\magnification=1000
\hsize=11.7cm
\vsize=18.9cm
\lineskip2pt \lineskiplimit2pt
\nopagenumbers

\hoffset=-1truein
\voffset=-1truein

\advance\voffset by 4truecm
\advance\hoffset by 4.5truecm

\newif\ifentete

\headline{\ifentete\ifodd	\count0 
      \rlap{\head}\hfill\tenrm\llap{\the\count0}\relax
    \else
        \tenrm\rlap{\the\count0}\hfill\llap{\head} \relax
    \fi\else
\global\entetetrue\fi}

\def\entete#1{\entetefalse\gdef\head{#1}} 
\entete{}

\input amssym.def
\input amssym.tex

\def\-{\hbox{-}}
\def\.{{\cdot}}
\def\O{{\cal O}}
\def\K{{\cal K}}
\def\F{{\cal F}}

\def\P{{\cal P}}

\def\G{{\cal G}}
\def\T{{\cal T}}

\def\I{{\cal I}}

\def\U{{\cal U}}
\def\H{{\cal H}}

\def\ch{\frak c\frak h}
\def\ad{\frak a\frak c}

\def\Gr{\frak G\frak r}

\def\Fct{\frak F\frak c\frak t}
\def\Nat{\frak N\frak a\frak t}

\def\int{\frak i\frak n\frak t}
\def\cat{\frak c\frak a\frak t}

\def\qq{\quad{\rm and}\quad}

\def\mod{\frak m\frak o\frak d}

\def\too{\longrightarrow}
\def\aut{\frak a\frak u\frak t}

\def\Set{\frak S\frak e\frak t}
\def\Loc{\frak L\frak o\frak c}
\def\loc{\frak l\frak o\frak c}

 3
 2
\font\large=cmr10  scaled \magstep 2
 2
\font\larti=cmti10  scaled \magstep 2
 1
 2

\font\cds=cmr7
\font\cdt=cmti7

\count0=1

\centerline{\large Categorizations of limits of  Grothendieck groups}
\medskip
\centerline{\large over a Frobenius {\larti P}-category}
\medskip
\centerline{\bf Lluis Puig }
\medskip
\noindent 
\centerline{\cds CNRS, Institut de Math\'ematiques de Jussieu, lluis.puig@imj-prg.fr}
\par
\noindent
\centerline{\cds 6 Av Bizet, 94340 Joinville-le-Pont, France}

\medskip
\noindent
{\bf Abstract:} {\cds  In~[9, Ch. 14] and [10]  we consider suitable inverse limits of Grothendieck groups of categories 
of modules in characteristics {\cdt p} and zero, obtained from a {\cdt folded Frobenius P-category\/} 
${\scriptstyle (\F,\widehat\aut_{\F^{^{\rm sc}}})}$ [10,~2.8], which covers the case of the {\cdt Frobenius P-categories\/} associated with {\cdt blocks\/};
moreover, in~[13] we show that a {\cdt folded Frobenius P-category\/}  is actually equivalent to  the choice of a 
{\cdt regular central ${\scriptstyle k^*\-}$extension \/}  ${\scriptstyle \hat\F^{^{\rm sc}}}$ of~${\scriptstyle \F^{^{\rm sc}}}\!$ [9,~11.2].  Here, taking advantage of the existence of a {\cdt perfect ${\scriptstyle \F^{^{\rm sc}}\!\-}$locality} 
${\scriptstyle \P^{^{\rm sc}}}$, recently proved in  [3], [5] and [11], we  exhibit those inverse limits as the  true 
{\cdt Grothendieck groups\/} of the categories of ${\scriptstyle \K_*\hat G\-}$ and ${\scriptstyle k_*\hat G\-}$modules for a suitable ${\scriptstyle k^*\-}$group  ${\scriptstyle \hat G}$ associated to  the ${\scriptstyle k^*\-}$category  ${\scriptstyle \hat\P^{^{\rm sc}}}\!$ 
  obtained from
${\scriptstyle \P^{^{\rm sc}}}\!$ and ${\scriptstyle \hat\F^{^{\rm sc}}}$. }

\bigskip
\noindent
{\bf £1. Introduction }
\medskip

£1.1. Let $p$ be a prime number and $\O$ a complete discrete valuation
ring with a {\it field of quotients\/} $\K$ of characteristic zero and a 
{\it residue field\/} $k$ of characteristic $p\,;$ we assume that $k$ is algebraically closed and that $\K$ contains ``enough'' roots of unity for the finite family of finite groups we will consider. Let $G$ be a finite group, $b$ a 
{\it block\/} of $G$ --- namely a primitive idempotent in the center $Z(\O G)$ of the group $\O\-$algebra --- and $(P,e)$ a maximal {\it Brauer $(b,G)\-$pair\/} [9,~1.16]; recall that the {\it Frobenius $P\-$category\/} $\F_{\!(b,G)}$ associated with~$b$ is the subcategory of the category of finite groups where the objects are all the subgroups of $P$ and, for any pair of subgroups $Q$ and $R$ of~$P\,,$ the morphisms $\varphi$ from $R$ to $Q$ are the group homomorphisms $\varphi\,\colon R\to Q$ induced by the conjugation of some element $x\in G$ fulfilling
$$(R,g)\i (Q,f)^x
\eqno £1.1.1\phantom{.}$$
where $(Q,f)$ and $(R,g)$ are the corresponding Brauer $(b,G)\-$pairs contained in $(P,e)$ [9, Ch.~3].

\medskip
£1.2. In~[9, Ch. 14] we consider a suitable inverse limit of Grothendieck groups of categories of modules in characteristic $p$ obtained from  $\F_{\!(b,G)}\,,$ which according to Alperin's Conjecture should be isomorphic to the 
 Grothendieck group of the category of finitely dimensional $kGb\-$modules.
In [10] we ge-neralize this construction in two directions. On the one hand,  we
also consider a suitable inverse limit of Grothendieck groups of categories of modules in characteristic zero  which again, according to Alperin's Conjecture, should be isomorphic to the  Grothendieck group of the category of finitely dimensional $\K Gb\-$modules. On the other hand, with the introduction of the 
{\it folded Frobenius $P\-$categories\/} [10, \S2], we are able to extend all these
constructions to any {\it folded Frobenius $P\-$category.\/}

\medskip
£1.3. Let us recall our definitions. Denoting by  $P$ a finite $p\-$group, by  $\frak i\Gr$ the category formed by the finite groups and  the injective group  homomorphisms, and  by $\F_{\!P}$ the subcategory of  $\frak i\Gr$ where the objects are all the  subgroups of $P$ and the morphisms are the group homomorphisms induced by the conjugation by elements of $P\,,$  a {\it Frobenius  $P\-$category\/} $\F$ is a subcategory 
of $\frak i\Gr$ containing $\F_{\!P}$ where the objects are all the  subgroups of $P$
and the morphisms fulfill the following three conditions [9, 2.8 and Proposition~2.11]
\smallskip
\noindent
£1.3.1\quad {\it For any subgroup $Q$ of $P$ the inclusion functor $(\F)_Q\to 
(\frak i\Gr)_Q$ is full.\/}
\smallskip
\noindent
£1.3.2\quad {\it $\F_P (P)$ is a Sylow $p\-$subgroup of $\F (P)\,.$\/}
\smallskip
\noindent
£1.3.3\quad {\it For any  subgroup $Q$ of $P$ such that we have $\xi \big(C_P (Q)\big)
= C_P\big(\xi (Q)\big)$ whenever $\xi\,\colon Q\.C_P(Q)\to P$ is an $\F\-$morphism,
any $\F\-$morphism $\varphi\,\colon Q\to P$ and any subgroup $R$ of $N_P\big(\varphi(Q)\big)$ containing $\varphi (Q)$ such that $\F_P(Q)$ contains the action of $\F_R \big(\varphi(Q)\big)$ over $Q$ via $\varphi\,,$ there is an $\F\-$morphism 
$\zeta\,\colon R\to P$ fulfilling $\zeta\big(\varphi (u)\big) = u$ for any $u\in Q\,.$\/}

\medskip
£1.4. Moreover, we say that a subgroup $Q$ of $P$ is {\it $\F\-$selfcentralizing\/} if we have
$$C_P\big(\varphi (Q))\i \varphi (Q)
\eqno £1.4.1\phantom{.}$$
 for any $\varphi \in \F (P,Q)\,,$ and we denote by $\F^{^{\rm sc}}$ the {\it full\/} subcategory of $\F$ over the set of  $\F\-$selfcentralizing subgroups of $P\,.$ We call 
{\it $\F^{^{\rm sc}}\!\-$chain\/} any functor 
 $\frak q\,\colon \Delta_n\to \F^{^{\rm sc}}$ where the $n\-$simplex $\Delta_n$ is considered as a category with the morphisms defined by the order [9,~A2.2]; we denote by $\Fct (\Delta_n, \F^{^{\rm sc}})$ this set of functors and by 
 $\ch^*(\F^{^{\rm sc}})$ the category  where the objects are all the 
$\F^{^{\rm sc}}\!\-$chains $(\frak q,\Delta_n)\,,$ where $n$ runs on $\Bbb N\,,$ and the morphisms from 
$\frak q\,\colon \Delta_n\to \F^{^{\rm sc}}$ to another $\F^{^{\rm sc}}\!\-$chain  $\frak r\,\colon \Delta_m\to \F^{^{\rm sc}}$ are the pairs $(\nu,\delta)$ 
formed by an {\it order preserving map\/} $\delta\,\colon \Delta_m\to \Delta_n$ and by a natural isomorphism 
$\nu\,\colon \frak q\circ\delta\cong \frak r$  [9,~A2.8]. Recall that we have a canonical functor
$$\aut_{\F^{^{\rm sc}}} : \ch^*(\F^{^{\rm sc}})\too \Gr
\eqno £1.4.2$$
mapping  any $\F^{^{\rm sc}}\!\-$chain $\frak q\,\colon \Delta_n\to \F^{^{\rm sc}}$ to the group of natural automorphisms 
of $\frak q$ [9, Proposition~A2.10].

\medskip
£1.5. Recall that a {\it $k^*\-$group\/} $\hat G$ is a group endowed with an
injective group homomorphism $\theta\,\colon k^*\to Z(\hat G)$~[7,~\S5],
 that $G = \hat G/\theta (k^*)$ is the {\it $k^*\-$quotient\/} of~$\hat G$ 
 and that a $k^*\-$group homomorphism is a group
homomorphism which preserves the multiplication by $k^*\,;$ let us denote by $k^*\-\Gr$ the category of $k^*\-$groups with finite $k^*\-$quotient.
Then, a {\it folded Frobenius $P\-$category\/} $(\F,\widehat\aut_{\F^{^{\rm sc}}})$ is a pair formed by a Frobenius $P\-$category $\F$ and,  by  a functor 
$$\widehat\aut_{\F^{^{\rm sc}}} : \ch^*(\F^{^{\rm sc}})\too k^*\-\Gr
\eqno £1.5.1\phantom{.}$$
lifting the canonical functor  $\aut_{\F^{^{\rm sc}}}\,;$ in the case of the {\it Frobenius $P\-$category\/} $\F_{\!(b,G)}$ above, we already know that the situation provides $k^*\-$groups $\hat\F_{\!(b,G)}(Q)$ lifting $\F_{\!(b,G)}(Q)$ for any 
$\F\-$selfcentralizing subgroup $Q$ of $P$ [9,~7.4] and in [9, Theorem~11.32] we prove the existence of a  lifting 
$\widehat\aut_{\F_{\!(b,G)}^{^{\rm sc}}}$ of $\aut_{\F_{\!(b,G)}^{^{\rm sc}}}$ extending them; note that  in 
[9, Theorem~11.32] we may assume that $k$ is just the closure of the prime subfield. But in [13,~Theorem~3.7] we prove  that any {\it folder structure\/} on $\F$ comes from  an essentially unique {\it regular central $k^*\-$extension\/} 
$\hat\F^{^{\rm sc}}$ of~$\F^{^{\rm sc}}$ and, from now on, a  {\it folder structure\/} on $\F$ means
a {\it regular central $k^*\-$extension\/} $\hat\F^{^{\rm sc}}$ of~$\F^{^{\rm sc}}$; if this {\it regular central 
$k^*\-$extension\/} comes from a {\it regular central $\bar k^*\-$extension\/} where $\bar k$ is the algebraic closure
 of the prime subfield then we say that the {\it folder structure\/} is {\it finite\/}.

\medskip
\noindent
{\bf Lemme~£1.6.} {\it For any finite folder structure $\hat\F^{^{\rm sc}}$ of~$\F^{^{\rm sc}}$ there are a finite subfield
$\hat k$ of $k$ and a regular central $\hat k^*\-$extension 
$\skew4\hat{\hat\F}^{^{\rm sc}}$ of~$\F^{^{\rm sc}}$ such that
the extension of $\skew4\hat{\hat\F}^{^{\rm sc}}$ from $\hat k^*$
to $k^*$ is equivalent to $\hat\F^{^{\rm sc}}\!$.\/}

\medskip
\noindent
{\bf Proof:} We actually may assume that $k$ is the algebraic closure of the prime subfield; then, for any
$\F^{^{\rm sc}}\!\-$morphism $\varphi\,\colon R\to Q\,,$ choose a lifting
$\hat\varphi\,\colon R\to Q$ of $\varphi$ in $\hat\F^{^{\rm sc}}\! (Q,R)\,;$
thus, for any pair of $\F^{^{\rm sc}}\!\-$morphisms $\varphi\,\colon R\to Q$
and $\psi\,\colon T\to R$ we get
$$\hat\varphi \circ\hat\psi = \lambda_{\varphi,\psi}\.\widehat{\varphi\circ\psi}
\eqno £1.6.1\phantom{.}$$
for a suitable {\it finite\/} family $\{\lambda_{\varphi,\psi}\}_{\varphi,\psi}$
of elements of $k^*$ which are algebraic; hence, the subfield~$\hat k$ of $k$
generated by this family is finite and it is clear that we can define a regular central 
$\hat k^*\-$extension $\skew4\hat{\hat\F}^{^{\rm sc}}$ of~$\F^{^{\rm sc}}$
contained in $\hat\F^{^{\rm sc}}\! $ by setting 
$$\skew4\hat{\hat\F}^{^{\rm sc}} (Q,R) = \bigcup_{\varphi\in \F^{^{\rm sc}}\! (Q,R)} \hat k^*\.\hat\varphi \i \hat\F^{^{\rm sc}}\! (Q,R)
\eqno £1.6.2\phantom{.}$$
and that this inclusion induces a bijection
$$k^*\times_{\hat k^*} \skew4\hat{\hat\F}^{^{\rm sc}} (Q,R) \cong  
\hat\F^{^{\rm sc}}\! (Q,R)
\eqno £1.6.3.$$

\medskip
£1.7. Note that a {\it regular central $k^*\-$extension\/} $\hat\F^{^{\rm sc}}$ 
of~$\F^{^{\rm sc}}$ induces a {\it regular central $k^*\-$extension\/} 
$\skew4\hat{\widetilde\F}^{^{\rm sc}}$ of the {\it exterior quotient\/} 
$\widetilde\F^{^{\rm sc}}$ of~$\F^{^{\rm sc}}$ --- namely of the quotient 
of~$\F^{^{\rm sc}}$ by the inner automorphisms of the objects [9,~6.1]; 
moreover, assuming that the {\it folder structure\/} $\hat\F^{^{\rm sc}}$ is {\it finite\/},
choosing a finite subfield $\hat k$ of $k$ and a  regular central $\hat k^*\-$extension $\skew4\hat{\hat\F}^{^{\rm sc}}$ 
of~$\F^{^{\rm sc}}$ as above, and denoting by $\skew4\hat{\skew4\hat{\widetilde\F}}^{^{\rm sc}}$ regular 
$\hat k^*\-$extension of the {\it exterior quotient\/} $\widetilde\F^{^{\rm sc}}\!$, if follows easily from [13,~Proposition~3.5] that the inclusion $\widetilde\F^{^{\rm sc}}_{\!P}\i \widetilde\F^{^{\rm sc}}$ can be lifted to a {\it faithful\/} functor 
$\widetilde\F^{^{\rm sc}}_{\!P}\to \skew4\hat{\skew4\hat{\widetilde\F}}^{^{\rm sc}}$ and then, with the terminology introduced in [12,~2.2], $\skew4\hat{\skew4\hat{\widetilde\F}}^{^{\rm sc}}$ becomes a 
{\it $\widetilde\F^{^{\rm sc}}_{\!P}\-$category\/} and fulfills the finiteness condition
in~[12,~4.1]; at this point, it follows from [12,~Proposition~4.6] that $\skew4\hat{\skew4\hat{\widetilde\F}}^{^{\rm sc}}$ becomes a {\it multiplicative $\widetilde\F^{^{\rm sc}}_{\!P}\-$category\/} since it inherits from $\widetilde\F^{^{\rm sc}}$ both conditions in this proposition.
\eject

\medskip
£1.8. On the other hand,  in [3], [5] and [11]  it has been recently proved that there exists a unique {\it perfect $\F^{^{\rm sc}}\!\-$locality\/} $\P^{^{\rm sc}}$ [9,~17.4 and~17.13]; more explicitly, denote by $\T^{^{\rm sc}}_P$ the category  where the objects are all the {\it $\F\-$self-centralizing\/} subgroups of $P\,,$ where the set of morphisms from $R$ to $Q$ is the {\it $P\-$transporter\/} $T_P (R,Q)$  for a pair of $\F\-$selfcentralizing subgroups $Q$ and $R$ of~$P\,,$ and where the composition is induced by the product in $P\,;$ then, there is a unique {\it Abelian extension\/} 
$\pi^{_{\rm sc}}\,\colon\P^{^{\rm sc}}\to \F^{^{\rm sc}}$ of  $\F^{^{\rm sc}}$ endowed with a {\it faithful\/} functor 
$\tau^{_{\rm sc}}\,\colon \T^{^{\rm sc}}_P \to \P^{^{\rm sc}}$  in such a way that the composition 
$\pi^{_{\rm sc}}\circ \tau^{_{\rm sc}}\,\colon \T^{^{\rm sc}}_P\to \F^{^{\rm sc}}$ is the canonical functor defined by the conjugation in~$P\,,$  that~$\P^{^{\rm sc}}\! (Q)$ endowed with 
$\tau^{_{\rm sc}}_Q\,\colon N_P(Q)\to \P^{^{\rm sc}}\! (Q)$ and 
$\pi^{_{\rm sc}}_Q \,\colon \P^{^{\rm sc}}\! (Q)\to \F^{^{\rm sc}} (Q)$ is an {\it $\F\-$localizer\/} of $Q$ for any
$\F\-$selfcentralizing subgroup $Q$ of $P$ fully normalized in $\F$ [9,~Theorem~18.6], and that $Z(R)$ acts {\it regularly\/} over the {\it fibers\/} of the map
$$\P^{^{\rm sc}}(Q,R)\too \F^{^{\rm sc}}(Q,R)
\eqno £1.8.1\phantom{.}$$
induced by $\pi^{_{\rm sc}}$ [9,~17.7] for any pair of $\F\-$selfcentralizing subgroups $Q$ and~$R$ of~$P\,.$

\medskip
£1.9.  Then, the so-called {\it $\F\-$localizing  functor\/} considered in [10, 3.2.1] 
$$\loc_{\F^{^{\rm sc}}} : \ch^*(\F^{^{\rm sc}})\too \widetilde{\Loc}
\eqno £1.9.1\phantom{.}$$
is  actually just a {\it quotient\/} of the canonical functor
$$\aut_{\P^{^{\rm sc}}} : \ch^* (\P^{^{\rm sc}})\too \Gr
\eqno £1.9.2\phantom{.}$$
mapping  any $\P^{^{\rm sc}}\!\-$chain $\frak q\,\colon \Delta_n\to \P^{^{\rm sc}}$ to the group of natural automorphisms 
of $\frak q$ [9, Proposition~A2.10]; moreover, a {\it regular central $k^*\-$extension \/} $\hat\F^{^{\rm sc}}$ 
of $\F^{^{\rm sc}}$  determines {\it via\/}~$\pi^{_{\rm sc}}$ a {\it regular central $k^*\-$extension \/} $\hat\P^{^{\rm sc}}$ of $\P^{^{\rm sc}}$ and, once again,
the {\it faithful\/} functor $\tau^{_{\rm sc}}\,\colon \T^{^{\rm sc}}_P \to \P^{^{\rm sc}}$ can be lifted to a {\it faithful\/} functor $\hat\tau^{_{\rm sc}}\,\colon \T^{^{\rm sc}}_P \to \hat\P^{^{\rm sc}}$ [13,~Proposition~3.5]; hence,  the corresponding functor
$$\widehat{\loc}_{\F^{^{\rm sc}}} : \ch^* (\F^{^{\rm sc}})\too k^*\-\widetilde\Loc
\eqno £1.9.3\phantom{.}$$
considered in [10, 3.3.1] is  presently just a {\it quotient\/} of the obvious canonical functor
$$\aut_{\hat\P^{^{\rm sc}}} : \ch^* (\hat\P^{^{\rm sc}})\too k^*\-\Gr
\eqno £1.9.4\phantom{.}$$
mapping  any $\hat\P^{^{\rm sc}}\!\-$chain $\hat\frak q\,\colon \Delta_n\to \hat\P^{^{\rm sc}}$ to the $k^*\-$group of natural automorphisms of $\hat\frak q$ [9, Proposition~A2.10], and  we simply set $\aut_{\hat\P^{^{\rm sc}}} (\hat\frak q,\Delta_n) = \hat\P^{^{\rm sc}}\! (\hat\frak q)\,.$

\medskip
£1.10. In this situation, considering the {\it contravariant\/} functors
$$\frak g_\K : k^*\-{\Gr}\too \O\-\mod\qq \frak g_k : k^*\-{\Gr}\too \O\-\mod
\eqno  £1.10.1\phantom{.}$$
mapping any $k^*\-$group  $\hat G$ with finite $k^*\-$quotient on the extensions to $\O$
$$\G_\K (\hat G) = \O\otimes_\Bbb Z \G^\Bbb Z_\K (\hat G)\qq \G_k (\hat G) = \O\otimes_\Bbb Z \G^\Bbb Z_k (\hat G)
\eqno £1.10.2\phantom{.}$$
\eject
\noindent
of the respective {\it  Grothendieck groups\/} $\G^\Bbb Z_\K (\hat G)$ and $\G^\Bbb Z_k (\hat G)$ of the categories of finitely dimensional $\K_*\hat G\-$ and $k_*\hat G\-$modules, and any $k^*\-$group homomorphism 
$\hat\theta\, \,\colon\hat G\to\hat G'$ on the corresponding {\it restriction\/} maps, with the the notation in~[10,~3.5] it is quite clear that 
$$\eqalign{\G_\K (\F,\widehat\aut_{\F^{^{\rm sc}}}) &= 
\lim_{\longleftarrow}\,(\frak g_\K\circ \widehat\loc_{\F^{^{\rm sc}}})\cong
\lim_{\longleftarrow}\,(\frak g_\K\circ \aut_{\hat\P^{^{\rm sc}}})\cr
\G_k (\F,\widehat\aut_{\F^{^{\rm sc}}}) &= 
\lim_{\longleftarrow}\,(\frak g_k\circ \widehat\loc_{\F^{^{\rm sc}}})\cong
\lim_{\longleftarrow}\,(\frak g_k\circ \aut_{\hat\P^{^{\rm sc}}})\cr}
\eqno £1.10.3.$$
Our purpose here is to exhibit these $\O\-$modules as the $\O\-$extensions of the very {\it Grothendieck groups\/} of suitable categories; as a matter of fact, we  find a $k^*\-$group $\hat G (\hat\P^{^{\rm sc}})\,,$ with finite $k^*\-$quotient, such that
$$\G_\K (\F,\widehat\aut_{\F^{^{\rm sc}}})\cong \G_\K \big(\hat G (\hat\P^{^{\rm sc}})\big)\qq
\G_k (\F,\widehat\aut_{\F^{^{\rm sc}}})\cong \G_k \big(\hat G (\hat\P^{^{\rm sc}})\big)
\eqno £1.10.4\phantom{.}$$
in a compatible way with the corresponding
 {\it decomposition maps\/} [10,~3.5.3]. We borrow our notation from [9] and [10].   
 
 \medskip
 \noindent
 {\bf Remark~£1.11.} Denoting by $\frak g_\K^{\Bbb Z^+}$ the {\it contravariant\/} functor mapping any $k^*\-$group 
 $\hat G$ with finite $k^*\-$quotient on the {\it positive\/} part $\G_\K^{\Bbb Z^+}\! (\hat G)$ of the very {\it Grothendieck group\/} of the category of finitely dimensional 
 $\K_*\hat G\-$modules  --- the part formed by the classes of the  $\K_*\hat G\-$modules  --- note that the {\it inverse limit\/} of the functor $\frak g^{\Bbb Z^+}_\K\circ \aut_{\hat\P^{^{\rm sc}}}$ still makes sense and actually it genera-tes 
 ${\displaystyle \lim_{\longleftarrow}}\,(\frak g^\Bbb Z_\K\circ \aut_{\hat\P^{^{\rm sc}}})\,.$ Indeed, denoting by~$R_{\hat\frak q}$ the class in $\G^{\Bbb Z^+}_\K \big(\hat\P^{^{\rm sc}}\! (\hat\frak q)\big)$ of the {\it regular\/} 
 $\K_*\hat\P^{^{\rm sc}}\! (\hat\frak q)\-$module $\K_*\hat\P^{^{\rm sc}}\! (\hat\frak q)$ for any 
 {\it $\hat\P^{^{\rm sc}}\-$chain\/}  $\hat\frak q\,\colon \Delta_n\to \hat\P^{^{\rm sc}}\,,$ and choosing
a multiple $m$ of all the orders $\vert\P^{^{\rm sc}}\! (\hat\frak q)\vert$ where $\hat\frak q$ runs over the set of  
{\it $\hat\P^{^{\rm sc}}\-$chains,\/} it is easily checked that  the family 
$R = \Big\{\displaystyle{m\over \vert\P^{^{\rm sc}}\! (\hat\frak q)\vert}\. R_{\hat\frak q}\Big\}_{\hat\frak q}$ belongs to 
${\displaystyle \lim_{\longleftarrow}}\, (\frak g^{\Bbb Z^+}_\K\circ \aut_{\hat\P^{^{\rm sc}}})$ and that, for any element $X$ of ${\displaystyle \lim_{\longleftarrow}}\,(\frak g^\Bbb Z_\K\circ \aut_{\hat\P^{^{\rm sc}}})\,,$  the sum of $X$ with a suitable 
multiple of~$R$  belongs  to ${\displaystyle \lim_{\longleftarrow}} (\frak g^{\Bbb Z^+}_\K\circ \aut_{\hat\P^{^{\rm sc}}})\,.$
The analogous argument also holds for~$\frak g^\Bbb Z_k\,.$

\bigskip
\noindent
{\bf £2. Categorization for the characteristic zero case}
\medskip
£2.1. Let $P$ be a finite $p\-$group and  $(\F,\hat\F^{^{\rm sc}})$  a 
{\it finite folded Frobenius $P\-$category\/};  denote by $\P$ and $\P^{^{\rm sc}}$ the respective {\it perfect\/} $\F\-$ and 
$\F^{^{\rm sc}}\-${\it localities\/} [11,~\S6~and~\S7] and by $\pi\,\colon \P\to \F$ and $\tau\,\colon \T_P\to \P$ the {\it structural functors} [9,~17.3]; then,  the {\it regular central $k^*\-$extension \/} $\hat\F^{^{\rm sc}}\!$ of~$\F^{^{\rm sc}}$  determines {\it via\/}~$\pi^{_{\rm sc}}\!$ a {\it regular central $k^*\-$extension \/} $\hat\P^{^{\rm sc}}\!$ of $\P^{^{\rm sc}}\!$ and we set (cf.~£1.10.3)
$$\G_\K (\hat\P^{^{\rm sc}}) = 
\lim_{\longleftarrow}\,(\frak g_\K\circ \aut_{\hat\P^{^{\rm sc}}})
\eqno £2.1.1;$$
that is to say, $\G_\K (\hat\P^{^{\rm sc}}) $ is the subset of elements
$$\{X_{(\hat\frak q,\Delta_n)}\}_{(\hat\frak q,\Delta_n)}\in \prod_{n\in \Bbb N}\,
\prod_{\hat\frak q\in \Fct (\Delta_n,\hat\P^{^{\rm sc}})} \G_\K \big(\hat\P^{^{\rm sc}}\! (\hat\frak q)\big)
\eqno £2.1.2\phantom{.}$$
such that, for any $\hat\P^{^{\rm sc}}\-$morphism $(\nu,\delta)\,\colon 
(\hat\frak q,\Delta_n)\to (\hat\frak r,\Delta_m)\,,$ they fulfill
$${\rm res}_{\aut_{\hat\P^{^{\rm sc}}}(\nu,\delta)}(X_{(\hat\frak r,\Delta_m)})
= X_{(\hat\frak q,\Delta_n)}
\eqno £2.1.3;$$
moreover, we assume that, for any extension $\K'$ of $\K\,,$ the scalar extension from $\K$ to $\K'$ induces an isomorphism between 
$\G_\K \big(\hat\P^{^{\rm sc}}\! (\hat\frak q)\big)$ and  
$\G_{\K'} \big(\hat\P^{^{\rm sc}}\! (\hat\frak q)\big)$ for any $n\in\Bbb N$
and any $\hat\frak q\in \Fct (\Delta_n,\hat\P^{^{\rm sc}})\,.$ Our purpose in this section is both to exhibit $\G_\K (\hat\P^{^{\rm sc}})$ as the extension to $\O$ of the very {\it Grothendieck group\/} of a  suitable subcategory of $\K\-\mod\-$valued 
{\it contravariant $k^*\-$functors\/} over the category $\hat\P^{^{\rm sc}}\,,$ and to show that this subcategory is equivalent
to the category of $\K_*\hat G\-$modules for a suitable $k^*\-$group $\hat G$ with finite $k^*\-$quotient.

\medskip
£2.2. First of all, it follows from [10,~Corollary~8.4] and from~£1.10.3 above that
$${\rm rank}_{\O}\big(\G_\K (\hat\P^{^{\rm sc}})\big) = 
\sum_{(\hat\frak q,\Delta_n)} (-1)^{n}\,{\rm rank}_{\O} 
\Big(\G_\K \big(\hat\P^{^{\rm sc}}\! (\hat\frak q)\big)\Big)
\eqno £2.2.1\phantom{.}$$
where $(\hat\frak q,\Delta_n)$ runs over a set of representatives for the set of isomorphism classes of $\ch^*(\hat\P^{^{\rm sc}})\-$objects such that 
$\hat\frak q (i-1,i)$ is {\it not\/} an isomorphism for any $1\le i\le n$ --- called
{\it regular $\ch^*(\hat\P^{^{\rm sc}})\-$objects\/} [9,~A5.2]. 
This formula suggests to consider the following {\it scalar product\/} in 
$\G_\K (\hat\P^{^{\rm sc}})\,;$ with the notation in~£2.1.2
 above, if $X = \{X_{(\hat\frak q,\Delta_n)}\}_{(\hat\frak q,\Delta_n)}$ and 
 $X' = \{X'_{(\hat\frak q,\Delta_n)}\}_{(\hat\frak q,\Delta_n)}$ are two elements of~$\G_\K (\hat\P^{^{\rm sc}})$ then we define
 $$\langle X,X'\rangle = \sum_{(\hat\frak q,\Delta_n)} (-1)^{n}\,
 \langle X_{(\hat\frak q,\Delta_n)}, X'_{(\hat\frak q,\Delta_n)} \rangle 
 \eqno £2.2.2\phantom{.}$$
 where $(\hat\frak q,\Delta_n)$ is running  over the same set and where, for such a 
 $(\hat\frak q,\Delta_n)\,,$ $\langle X_{(\hat\frak q,\Delta_n)}, X'_{(\hat\frak q,\Delta_n)} \rangle $ denotes the scalar product of $X_{(\hat\frak q,\Delta_n)}$ and $X'_{(\hat\frak q,\Delta_n)}$ in the Grothendieck group
 $\G_\K \big(\hat\P^{^{\rm sc}}\! (\hat\frak q)\big)\,.$ Note that there is a canonical bijection between  a set of representatives for the set of isomorphism classes of {\it regular $\ch^*(\hat\P^{^{\rm sc}})\-$objects\/} and a set of
representatives for the set of $\F\-$isomorphism classes of nonempty sets of
$\F\-$selfcentralizing subgroups of~$P$,  totally ordered by the inclusion.

 \medskip
 £2.3. Recall that the canonical group homomorphism $\O^*\to k^*$ admits
 a unique section $k^*\to \O^*\,;$ thus, the category $\K\-\mod$ admits
 an evident $k^*\-$action and a {\it contravariant $k^*\-$functor\/} 
 $\frak m\,\colon \hat\P^{^{\rm sc}}\!\to \K\-\mod$ is a functor such that   
  $\frak m[\lambda\.\hat x) = \lambda\.\frak m (\hat x)$
 for any $\hat\P^{^{\rm sc}}\!\-$morphism $\hat x\,\colon R\to Q$ and any 
 $\lambda\in k^*\,.$\break
 \eject
 \noindent
  Moreover, it is  clear that any {\it contravariant $k^*\-$functor\/} 
 $\frak m\,\colon \hat\P^{^{\rm sc}}\!\to \K\-\mod$ determines a new 
 {\it contravariant $k^*\-$functor\/}  [9,~A3.7.3]
 $$\frak m^{\ch} = \frak m\circ \frak v_{\hat\P^{^{\rm sc}}}: 
 \ch^* (\hat\P^{^{\rm sc}})\too \K\-\mod
 \eqno £2.3.1\phantom{.}$$
 sending any {\it $\hat\P^{^{\rm sc}}\!\-$chain\/} 
 $\hat\frak q\,\colon \Delta_n\to \hat\P^{^{\rm sc}}\!$ to 
 $\frak m\big(\hat\frak q (0)\big)$ and any 
 $\ch^*(\hat\P^{^{\rm sc}})\-$morphism
 $$(\hat x,\delta) : (\hat\frak r,\Delta_m)\too (\hat\frak q,\Delta_n)
 \eqno £2.3.2,$$
where $\delta\,\colon \Delta_n\to \Delta_m$ is an order-preserving map
and $\hat x\,\colon \hat\frak r\circ\delta\cong \hat\frak q$ a 
{\it natural isomorphism\/}, to the $\K\-$linear map 
$$\frak m \big(\hat x_0\circ \hat r (0\bullet \delta (0)\big)
:  \frak m\big(\hat\frak q (0)\big)\too \frak m\big(\hat\frak r (0)\big)
\eqno £2.3.3.$$

\medskip
£2.4. In particular, in the case where $n = m$ and $\delta$ is the identity map, we get a $k^*\-$compatible action  over 
$\frak m\big(\hat\frak q (0)\big)$ of the $k^*\-$group 
$\aut_{\hat\P^{^{\rm sc}}\!} (\hat\frak q) = \hat\P^{^{\rm sc}}\! (\hat\frak q)$; that is to say, 
$\frak m\big(\hat\frak q (0)\big)$ becomes a 
$\K_*\hat\P^{^{\rm sc}}\! (\hat\frak q)\-$module. Similarly, 
$\frak m\big(\hat\frak r (0)\big)$ becomes a $\K_*\hat\P^{^{\rm sc}}\! 
(\hat\frak r)\-$module and, {\it via\/} the $k^*\-$group homomorphism 
$$\aut_{\hat\P^{^{\rm sc}}\!} (\hat x,\delta) : \hat\P^{^{\rm sc}}\! 
(\hat\frak r)\too \hat\P^{^{\rm sc}}\! (\hat\frak q)
\eqno £2.4.1,$$
$\frak m\big(\hat\frak q (0)\big)$ also becomes the $\K_*\hat\P^{^{\rm sc}}\! 
(\hat\frak r)\-$module ${\rm Res}_{\aut_{\hat\P^{^{\rm sc}}\!} (\hat x,\delta)}
\Big(\frak m\big(\hat\frak q (0)\big)\Big)$ and then the  $\K\-$linear map~£2.3.3 
 is clearly a  $\K_*\hat\P^{^{\rm sc}}\! (\hat\frak r)\-$module homomorphism
 $$\frak m^{\ch} (\hat x,\delta) :  
 {\rm Res}_{\aut_{\hat\P^{^{\rm sc}}\!} (\hat x,\delta)}
\Big(\frak m\big(\hat\frak q (0)\big)\Big)\too \frak m\big(\hat\frak r (0)\big)
\eqno £2.4.2.$$

\medskip
£2.5. Similarly, any {\it natural map\/} $\mu\,\colon \frak m\to \frak m'$ between
 {\it contravariant $k^*\-$ functors\/} $\frak m$ and $\frak m'$ from 
 $\hat\P^{^{\rm sc}}\!$ to $\K\-\mod$ determines  a new 
 {\it natural map\/} 
 $$\mu^{\ch} = \mu *\frak v_{\hat\P^{^{\rm sc}}} : 
 \frak m^{\ch}\too \frak m'^{\ch}
 \eqno £2.5.1\phantom{.}$$
 sending any {\it $\hat\P^{^{\rm sc}}\!\-$chain\/} 
 $\hat\frak q\,\colon \Delta_n\to \hat\P^{^{\rm sc}}\!$ to the $\K\-$linear map
 $$\mu_{\hat\frak q (0)} : \frak m\big(\hat\frak q (0)\big)\too
 \frak m'\big(\hat\frak q (0)\big)
 \eqno £2.5.2;$$
then, it follows from the naturalness of $\mu$ that this map is actually
a $\K_*\hat\P^{^{\rm sc}}\! (\hat\frak q)\-$ module homomorphism.
Let us denote by $\Nat (\frak m',\frak m)$ the $\K\-$module of 
{\it natural maps\/} from $\frak m$ to $\frak m'\,.$

\medskip
£2.6. We are interested in the {\it contravariant $k^*\-$functors\/} 
 $\frak m\,\colon \hat\P^{^{\rm sc}}\!\to \K\-\mod$ --- called {\it reversible\/} --- mapping {\it any\/}  $\hat\P^{^{\rm sc}}\!\-$morphism $\hat x\,\colon R\to Q$ on a $\K\-$linear isomorphism 
 $$\frak m(\hat x) : \frak m(Q)\cong \frak m(R)
 \eqno £2.6.1\,;$$
  in this case, it is quite clear that the {\it contravariant $k^*\-$functor\/}  
 $\frak m^{\ch}$ also maps any $\ch^* (\hat\P^{^{\rm sc}})\-$morphism
 on a $\K\-$linear isomorphism. Conversely, 
 \smallskip
 \noindent
 £2.6.2\quad {\it A  contravariant $k^*\-$functor
 $\,\frak x : \ch^* (\hat\P^{^{\rm sc}})\too \K\-\mod$  mapping any
 $\ch^* (\hat\P^{^{\rm sc}})\-$morphism on a $\K\-$linear isomorphism comes from a reversible contravariant $k^*\-$functor.\/}
 \eject
 
 \medskip
 £2.7. Indeed, for any $\F\-$selfcentralizing subgroup $Q$ of $P\,,$ considering the  
{\it $\hat\P^{^{\rm sc}}\!\-$chain\/} $\hat\frak q_Q\,\colon\Delta_0\to 
\hat\P^{^{\rm sc}}\!$ sending $0$ to $Q\,,$ we set $\frak m (Q) = 
\frak x (\hat\frak q_Q,\Delta_0)\,;$ moreover, for any 
$\hat\P^{^{\rm sc}}\-$morphism $\hat x\,\colon R\to Q\,,$  considering the 
{\it $\hat\P^{^{\rm sc}}\!\-$chain\/} 
$\hat\frak q_{\hat x}\,\colon\Delta_1\to \hat\P^{^{\rm sc}}\!$ mapping $0$ on $R\,,$ $1$ on $Q$ and the $\Delta_1\-$morphism $0\bullet 1$ on $\hat x\,,$
we have the evident  $\ch^*(\hat\P^{^{\rm sc}})\-$morphisms
 $$\eqalign{({\rm id}_R,\delta_1^0) : (\hat\frak q_{\hat x},\Delta_1)\too (\hat\frak q_R,\Delta_0)\cr 
 ({\rm id}_Q,\delta_0^0) : (\hat\frak q_{\hat x},\Delta_1)\too (\hat\frak q_Q,\Delta_0)\cr}
 \eqno £2.7.1\phantom{.}$$
 and  we set $\frak m (\hat x) = \frak x( ({\rm id}_R,\delta_0^0))^{-1}\circ 
\frak x ({\rm id}_Q,\delta_1^0)\,.$ This correspondence is actually a 
$k^*\-$functor since, for another $\hat\P^{^{\rm sc}}\-$morphism 
$\hat y\,\colon T\to R\,,$ considering the new {\it $\hat\P^{^{\rm sc}}\!\-$chain\/} $\hat\frak r\,\colon\Delta_2\to \hat\P^{^{\rm sc}}\!$ sending $0$ to $T\,,$ $1$ to $R\,,$ $2$ to $Q\,,$ $0\bullet 1$ to $\hat y$ and $1\bullet 2$ 
to~$\hat x\,,$ and extending the notation above, we get the evident
commutative diagram in the category $\ch^*(\hat\P^{^{\rm sc}})$
$$\matrix{\hat\frak q_T& &\hat\frak q_R& &\hat\frak q_Q\cr
\uparrow&\nwarrow\hskip-10pt\nearrow&\uparrow&\nwarrow\hskip-10pt\nearrow&\uparrow\cr
\hat\frak q_{\hat y}&&\hat\frak q_{\hat x\.\hat y}&& \hat\frak q_{\hat x}\cr
&\nwarrow&\uparrow&\nearrow\cr
&&\hat\frak r\cr}
\eqno £2.7.2\phantom{.}$$
which the functor $\frak x$ sends to a commutative diagram of isomorphisms 
in $\K\-\mod\,;$ then, it is easily checked that $\frak m (\hat x\.\hat y) = 
\frak m (\hat y)\circ \frak m (\hat x)$ and that $\frak m^{\ch} = \frak x\,.$

\medskip
£2.8. If $\frak m$ and $\frak m'$ are {\it reversible contravariant $k^*\-$functors\/} from the $k^*\-$ca-tegory 
$\hat\P^{^{\rm sc}}\!$ to~$\K\-\mod\,,$ a {\it natural map\/} 
$\mu\,\colon \frak m\to \frak m'$ is a correspondence sending any 
$\F\-$selfcentralizing subgroup $Q$ of $P$
to a $\K\-$linear map $\mu_Q\,\colon \frak m (Q)\to \frak m' (Q)$
in such a way that, for any $\hat\P^{^{\rm sc}}\!\-$morphism 
$\hat x\,\colon R\to Q\,,$ we have the commutative diagram
$$\matrix{\frak m(Q)&\buildrel \frak m (\hat x)\over\cong &\frak m(R)\cr
\hskip-20pt{\scriptstyle \mu_Q}\hskip5pt\big\downarrow&\phantom{\Big\downarrow}&\big\downarrow\hskip5pt{\scriptstyle \mu_R}\hskip-20pt\cr
\frak m' (Q)&\buildrel \frak m' (\hat x)\over\cong &\frak m' (R)\cr}
\eqno £2.8.1;$$
equivalently, considering the  {\it reversible contravariant functor\/}
$$\frak m^*\otimes_\K \frak m' : \P^{^{\rm sc}}\!\too \K\-\mod
\eqno £2.8.2\phantom{.}$$
mapping $Q$ on ${\rm Hom}_\K \big(\frak m(Q),\frak m'(Q)\big)$
and $\hat x\,\colon R\to Q$ on the $\K\-$linear map
$${\rm Hom}_\K \big(\frak m(Q),\frak m'(Q)\big)\too 
{\rm Hom}_\K \big(\frak m(R),\frak m'(R)\big)
\eqno £2.8.3\phantom{.}$$
sending any $\alpha\in {\rm Hom}_\K \big(\frak m(Q),\frak m'(Q)\big)$
to $\frak m' (\hat x)\circ \alpha\circ \frak m (\hat x)^{-1}\,,$ the commutativity of the diagrams above means that $\mu$ belongs to the inverse limit 
of~$\frak m^*\otimes_\K \frak m'\,;$ that is to say, we get
$$\Nat (\frak m',\frak m) = \lim_{\longleftarrow}\,
(\frak m^*\otimes_\K \frak m')
\eqno £2.8.4.$$

\medskip
£2.9. Moreover, let $\frak m\,\colon \hat\P^{^{\rm sc}}\!\to \K\-\mod$ be a
{\it reversible contravariant $k^*\-$ functor\/}; 
 with the notation in~£2.4 above, homomorphism~£2.4.2 becomes an isomorphism
 and therefore the restriction map induced by homomorphism~£2.4.1
 $$ \frak g_\K \big(\aut_{\hat\P^{^{\rm sc}}\!} (\hat x,\delta)\big) :
 \G_\K (\hat\P^{^{\rm sc}}\! (\hat\frak q)\big)\too \G_\K (\hat\P^{^{\rm sc}}\! (\hat\frak r)\big)
 \eqno £2.9.1\phantom{.}$$
  sends  the class $X_{\frak m (\hat\frak q (0))}$ in $\G_\K (\hat\P^{^{\rm sc}}\! (\hat\frak q)\big)$ of the 
 $\K_*\hat\P^{^{\rm sc}}\! (\hat\frak q)\-$module  $\frak m \big(\hat\frak q (0)\big)$ to the class 
 $X_{\frak m (\hat\frak r (0))}$  in $\G_\K (\hat\P^{^{\rm sc}}\! (\hat\frak r)\big)$ of the 
 $\K_*\hat\P^{^{\rm sc}}\! (\hat\frak r)\-$module 
 $\frak m \big(\hat\frak r (0)\big)\,.$ That is to say,\break
 the family
 $$\big\{X_{\frak m (\hat\frak q (0))}\big\}_{(\hat\frak q,\Delta_n)}
 \in \prod_{n\in\Bbb N}\,\prod_{\hat\frak q\in \Fct (\Delta_n,\hat\P^{^{\rm sc}})}
 \G_\K (\hat\P^{^{\rm sc}}\! (\hat\frak q)\big) 
 \eqno £2.9.2\phantom{.}$$
 fulfills condition~£2.1.3 and therefore it belongs to 
 $\G_\K (\hat\P^{^{\rm sc}})\,;$  let us denote by $X_\frak m$ this family which, clearly, only depends on the isomorphism class of~$\frak m\,.$

\medskip
£2.10.  Actually, any {\it reversible contravariant  $k^*\-$functor\/}  $\frak m\,\colon \hat\P^{^{\rm sc}}\!\to \K\-\mod$ is 
{\it naturally isomorphic\/} to a  {\it contravariant $k^*\-$functor\/} $\frak n$ from $\hat\P^{^{\rm sc}}\!$ to $\K\-\mod$ --- called {\it reduced\/} ---  mapping any $\F\-$selfcentralizing subgroup $Q$ of~$P$ on the {\it same\/} finite dimensional 
$\K\-$module $M$ and  any $\hat\P^{^{\rm sc}}\!\-$morphism $\hat\tau_{Q,R}(1)$ on~${\rm id}_M\,;$ indeed, setting 
$M = \frak m (P)\,,$ we define $\frak n$ sending any $\F\-$selfcentralizing subgroup $Q$ of~$P$ to $M$ and any
$\hat\P^{^{\rm sc}}\!\-$morphism $\hat x\,\colon R\to Q$ to the $\K\-$linear map
$$\frak m\big(\hat\tau_{P,R} (1)\big)^{-1}\circ \frak m (\hat x)\circ 
\frak m\big(\hat\tau_{P,Q} (1)\big) : M\cong M
\eqno £2.10.1$$
which, if $R$ is contained in $Q\,,$ clearly maps $\hat\tau_{Q,R}(1)$ on ${\rm id}_M\,;$ then,
we have the obvious {\it natural isomorphism\/} $\frak n\cong \frak m$ sending
$Q$ to 
$$\frak m\big(\hat\tau_{P,Q} (1)\big) : M\cong \frak m(Q)
\eqno £2.10.2.$$
Moreover, two of such {\it reduced contravariant  $k^*\-$functors\/} $\frak n$
and $\frak n'$ mapping $P$ on the same $\K\-$module $M$ are {\it naturally isomorphic\/} if and only if there is $s\in {\rm GL}_\K (M)$ fulfilling 
$\frak n' (\hat x) = \frak n (\hat x)^s$ for any $\hat\P^{^{\rm sc}}\!\-$morphism $\hat x\,\colon R\to Q\,.$

 \medskip
 £2.11. Finally, we call {\it $\K_*\hat\P^{^{\rm sc}}\!\-$module\/} any  {\it reversible contravariant  $k^*\-$fun-ctor\/}  $\frak m$ such that, setting $M = \frak m (P)\,,$ the $k^*\-$subgroup $\hat G(\frak m)$ of  ${\rm GL}_\K (M)$ generated 
 by~$\frak m\big(\hat\tau_{P,R} (1)\big)^{-1}\circ \frak m (\hat x)\circ \frak m\big(\hat\tau_{P,Q} (1)\big)\,,$ where  
 $\hat x\,\colon R\to Q$ runs over the set of  $\hat\P^{^{\rm sc}}\!\-$morphisms, has a finite $k^*\-$quotient $G(\frak m)\,;$ 
 it is clear that the {\it direct sum\/} of  {\it $\K_*\hat\P^{^{\rm sc}}\!\-$modules\/} 
 is a {\it $\K_*\hat\P^{^{\rm sc}}\!\-$module\/}; we denote by $\K_*\hat\P^{^{\rm sc}}\!\-\mod$
 the category formed by the {\it $\K_*\hat\P^{^{\rm sc}}\!\-$modules\/} and by the {\it natural maps\/} between them.

 \bigskip
 \noindent
 {\bf Proposition~£2.12.} {\it If $\,\frak m\,\colon \hat\P^{^{\rm sc}}\!\to \K\-\mod$ is a $\K_*\hat\P^{^{\rm sc}}\!\-$module then any  reversible contravariant  $k^*\-$functor $\frak n\,\colon \hat\P^{^{\rm sc}}\!\to  \K\-\mod$ contained in $\frak m$  is a $\K_*\hat\P^{^{\rm sc}}\!\-$module too and admits a complement in $\frak m\,.$ In particular, $\frak m$ is isomorphic 
 to a direct sum of simple $\K_*\hat\P^{^{\rm sc}}\!\-$modules.\/}
 \eject

 \medskip
 \noindent
 {\bf Proof:} We may assume that $\frak m$ and $\frak n$ are {\it reduced\/} and, setting $N = \frak n (P)$ and 
 $M = \frak m (P)\j N\,,$ it is clear that $\hat G (\frak m)\i {\rm GL}_\K (M)$  stabilizes $N$ and that its image in ${\rm GL}_\K (N)$ 
 contains~$\frak n (\hat x)$ where  $\hat x\,\colon R\to Q$ runs over the set of  $\hat\P^{^{\rm sc}}\!\-$morphisms, 
 so that $G (\frak n)$ is finite too. Moreover, since $M$ is a semisimple $\K_*\hat G (\frak m)\-$module and $N$
 is a $\K_*\hat G (\frak m)\-$submodule, we have $M\cong N\oplus N'$ for a suitable $\K_*\hat G (\frak m)\-$submodule $N'$
 and it is quite clear that $N'$ defines a $\K_*\hat\P^{^{\rm sc}}\!\-$submodule
 of $\frak m$ which is a complement of $\frak n\,.$

 \medskip
 £2.13.  Thus, in the category  of $\K_*\hat\P^{^{\rm sc}}\!\-$modules any object is a direct sum of simple objects
 and therefore it makes sense to talk about  the {\it Grothendieck group\/} of this category; we denote by 
 $\G (\K_*\hat\P^{^{\rm sc}}\!\-\mod)$ the extension to $\O$ of this Grothendieck group and then it is clear that the correspondence sending any {\it $\K_*\hat\P^{^{\rm sc}}\!\-$module\/} $\frak m$ to 
 $X_\frak m$ (cf.~£2.9) induces an $\O\-$module homomorphism
 $$\cat_\K : \G (\K_*\hat\P^{^{\rm sc}}\!\-\mod)\too \G_\K (\hat\P^{^{\rm sc}})
 \eqno £2.13.1.$$
Actually,  from the following result we will show that this $\O\-$module homomorphism  is injective proving, in particular,
that  the set of isomorphism classes of simple $\K_*\hat\P^{^{\rm sc}}\!\-$modules is finite.

 \medskip
 £2.14.     In this result we have to deal with {\it $n\-$cohomology groups\/} of the ca-tegory $\P^{^{\rm sc}}\!$ over  
 a {\it (reversible) contravariant functor\/}   $\frak m\,\colon \P^{^{\rm sc}}\!\to \K\-\mod\,;$ namely, for any 
 $n\in \Bbb N$ we set
 $$\Bbb C^n  (\P^{^{\rm sc}}\!,\frak m) = 
 \prod_{\frak q\in \Fct (\Delta_n,\P^{^{\rm sc}})} \frak m\big(\frak q (0)\big)
 \eqno £2.14.1\phantom{.}$$
 and consider the usual differential map $d^n\,\colon \Bbb C^n  (\P^{^{\rm sc}}\!,\frak m) \to 
 \Bbb C^{n+1}  (\P^{^{\rm sc}}\!,\frak m)$ [9, A3.11.2]; more generally, we denote by $\Bbb C^n_* (\P^{^{\rm sc}}\!,\frak m)$ the $\K\-$submodule of {\it stable\/} families, namely the families  $\{m_{\frak q}\}_{\frak q\in \Fct (\Delta_n,\P^{^{\rm sc}})}$
  such that $\big(\frak m (x_0)\big)(m_{\frak q'}) =  m_{\frak q}$ for any   {\it natural isomorphism\/}  $x\,\colon \frak q\cong \frak q'\,;$ the differential map restricts to a new differential map $d_*^n$ and we denote by 
  $\Bbb H_*^n (\P^{^{\rm sc}}\!,\frak m) $ the corresponding {\it stable $n\-$cohomology groups\/} [9,~A3.18]; similarly,
considering the {\it regular $\ch^*(\P^{^{\rm sc}})\-$ objects\/} in~£2.2 above, we denote by $\Bbb C_{\rm r}^n  (\P^{^{\rm sc}}\!,\frak m)$  
the quotient of $\Bbb C^n_* (\P^{^{\rm sc}}\!,\frak m)$ by all the elements vanishing over the {\it regular $\P^{^{\rm sc}}\-$chains\/} 
$\Fct_{\rm r} (\Delta_n,\P^{^{\rm sc}})$ [9,~A5.3]
and we get the {\it regular $n\-$cohomology groups\/} 
$\Bbb H_{\rm r}^n (\P^{^{\rm sc}}\!,\frak m) $ [9,~A5.6].

  \bigskip
 \noindent
 {\bf Theorem~£2.15.} {\it For any $\K\P^{^{\rm sc}}\!\-$module $\frak m$
 and any $n\ge 1$ we have 
 $$\Bbb H^n (\P^{^{\rm sc}}\!,\frak m) =  \Bbb H^n_* (\P^{^{\rm sc}}\!,\frak m) =  \Bbb H^n_{\rm r} (\P^{^{\rm sc}}\!,\frak m) = \{0\}
 \eqno £2.15.1.$$\/}
 
 \par
\noindent
{\bf Proof:}  The equality $\Bbb H^n (\P^{^{\rm sc}}\!,\frak m) = 
 \Bbb H^n_* (\P^{^{\rm sc}}\!,\frak m)$ follows from [9, Proposition~A4.13] and the equality 
 $\Bbb H^n_* (\hat\P^{^{\rm sc}}\!,\frak m) = 
 \Bbb H^n_{\rm r} (\hat\P^{^{\rm sc}}\!,\frak m)$ follows from [9, Proposition~A5.7]. More generally, denoting by $\I$ the  {\it interior structure\/} 
 of~$\P^{^{\rm sc}}\!$ [9,~1.3] mapping any $\F\-$selfcentralizing subgroup 
 $Q$ of $P$ on  $\tau^{_{\rm sc}}_{_Q} (Q)\,,$ it still follows from 
 [9, Proposition~A4.13] that $\Bbb H^n (\P^{^{\rm sc}}\!,\frak m)$ coincides with the {\it $\I\-$stable $n\-$cohomology group 
 $\Bbb H^n_{\I} (\P^{^{\rm sc}}\!,\frak m)$ of~$\P^{^{\rm sc}}\!$\/} 
 over $\frak m$ [9,~A3.18]; thus, it suffices to prove that, for  any~$n\ge 1\,,$  
 $\Bbb H^n_{\I }(\P^{^{\rm sc}}\!,\frak m) = \{0\}\,.$

 \smallskip
 We may assume that $\frak m$ is {\it reduced\/}; setting $M = \frak m (P)\,,$ since the $k^*\-$group $\hat G(\frak m)\i {\rm GL}_\K (M)$ has a finite $k^*\-$quotient,
 it  stabilizes an $\O\-$submodule $M^{_\O}$ of $M$ such that 
 $M\cong \K\otimes_\O M^{_\O}\,;$  then, $\hat G(\frak m)$ is contained
 in ${\rm GL}_\O (M^{_\O})$ and it suffices to define $\frak m^{_\O}(\hat x)
 = \frak m (\hat x)$  for any $\hat\P^{^{\rm sc}}\!\-$morphism $\hat x\,\colon R\to Q$ to get a {\it contravariant\/} functor $\frak m^{_\O}$  from $\P^{^{\rm sc}}$ to $\O\-\mod$ such that $\frak m\cong \K\otimes_\O \frak m^{_\O}\,;$ hence, it suffices to prove that, for  any $n\ge 1\,,$ we have
$$\Bbb H^n_{\I } (\P^{^{\rm sc}}\!,\frak m^{_\O}) = \{0\}
\eqno £2.15.2.$$
  
\smallskip
More precisely, setting  $\bar\frak m^{_\O} = k\otimes_\O \frak m^{_\O}$ and $\bar M^{^\O} = k\otimes_\O M^{^\O}\,,$ 
it suffices to prove that 
$$\Bbb H^n_{\I } (\P^{^{\rm sc}}\!,\bar\frak m^{_\O}) = \{0\}
\eqno £2.15.3;$$ 
indeed, as in~£2.14 above, for any $n\in \Bbb N$ we set
 $$\Bbb C^n  (\P^{^{\rm sc}}\!,\frak m^{_\O}) = 
 \prod_{\frak q\in \Fct (\Delta_n,\P^{^{\rm sc}})} 
 M^{^\O}
 \eqno £2.15.4\phantom{.}$$
 and denote by $\Bbb C^n_{\I }  (\P^{^{\rm sc}}\!,\frak m^{_\O})$ the $\O\-$submodule of {\it $\I\-$stable\/} families;  
 then, if equality~£2.15.3 holds and $c_0\in \Bbb C^n_{\I }  (\P^{^{\rm sc}}\!,\frak m^{_\O})$ is an
{\it  $n\-$cocycle\/}, denoting by $\varpi$ a generator of $J(\O)\,,$ we already have  that
$$c_0 \equiv d_{n-\!1}(a_0)\pmod \varpi
\eqno £2.15.5\phantom{.}$$ 
for a suitable  $a_0\in \Bbb C^{n-1}_{\I }  (\P^{^{\rm sc}}\!,\frak m^{_\O})\,,$ so that we have 
$c_0 - d^{n-\!1}_{\I }(a_0) = \varpi\.c_1$ for a unique $c_1\in \Bbb C^n_{\I }  (\P^{^{\rm sc}}\!,\frak m^{_\O})$ which is
again an {\it  $n\-$cocycle\/} since $M^{^\O}$ is a free $\O\-$module; thus, for any $i\in \Bbb N$ we inductively can define  $c_i\in \Bbb C^n_{\I }  (\P^{^{\rm sc}}\!,\frak m^{_\O})$ and 
$a_i\in \Bbb C^{n-1}_{\I }  (\P^{^{\rm sc}}\!,\frak m^{_\O})$ fulfilling 
$$c_i\equiv d^{n-\!1}_{\I } (a_i)\pmod \varpi\qq c_i - d^{n-\!1}_{\I } (a_i) 
= \varpi\. c_{i+\!1}
\eqno £2.15.6\phantom{.}$$
and then, according to the completeness of $\O\,,$ it is quite clear that
$$c_0 = d^{n-\!1}_{\I } \big(\sum_{i\in \Bbb N} \varpi^i\.a_i\big)
\eqno £2.15.7.$$

\smallskip
Now, denoting by $\frak h_0(\bar\frak m^{_\O})\,\colon 
\P^{^{\rm sc}}\to k\-\mod$ the {\it contravariant\/} $k^*\-$subfunctor of
 $\bar\frak m^{_\O}$ sending any $\F\-$selfcentralizing subgroup $Q$ of $P$ to  the fixed points of $Q$ in $\bar\frak m^{_\O}(Q)$ [9,~14.21], it is quite clear that the inclusion of $\frak h_0(\bar\frak m^{_\O})$ in $\bar\frak m^{_\O}$ induces  
 $\O\-$module isomorphisms
$$\Bbb C^n_{\I } \big(\P^{^{\rm sc}}\!,\frak h_0(\bar\frak m^{_\O})\big)\cong \Bbb C^n_{\I }(\P^{^{\rm sc}}\!,\bar\frak m^{_\O})
\eqno £2.15.8\phantom{.}$$
for any $n\in \Bbb N$ which are compatible with the differential maps; thus,
 it suffices to prove that $\Bbb H_{\I }^n \big(\P^{^{\rm sc}}\!,
 \frak h_0(\bar\frak m^{_\O})\big) = \{0\}$ for any $n\ge 1\,.$ 
 \eject

 \smallskip
 Moreover, the {\it $\I\-$exterior quotient\/} of $\P^{^{\rm sc}}\!$  [9,~1.3] coincides with the  {\it exterior quotient\/} 
 $\widetilde\F^{^{\rm sc}}\!$ (cf.~£1.7)  and then  the  {\it contravariant\/} functor $\frak h_0(\bar\frak m^{_\O})$ factorizes through the canonical functor~$\P^{^{\rm sc}}\to \widetilde\F^{^{\rm sc}}\,.$
On the other hand, as in~£1.7 above,  with the terminology introduced in [12,~2.2] $\P^{^{\rm sc}}\!$ becomes a {\it $\T^{^{\rm sc}}_{\!P}\-$category\/} and it fulfills the {\it finiteness condition\/} in~[12,~4.1]; at this point, it follows from [12,~Proposition~4.6] that $\P^{^{\rm sc}}\!$ is a 
{\it multiplicative $\T^{^{\rm sc}}_{\!P}\-$category\/} since, by [9, Propositions~24.2 and~24.4],  $\P^{^{\rm sc}}$ fulfills both conditions in this proposition. 

\smallskip
That is to say, the {\it additive cover\/} $\ad(\P^{^{\rm sc}})$ admits  {\it direct products\/} [12,~4.2];  in particular, we have a functor
$$\int_P : \P^{^{\rm sc}}\!\too \ad (\P^{^{\rm sc}})
\eqno £2.15.9\phantom{.}$$
sending any $\F\-$selfcentralizing subgroup $Q$ of $P$ to the {\it direct product\/} $Q\hat\times P$ 
in~$\ad (\P^{^{\rm sc}})\,;$ hence, setting $ i_{Q}^P = \tau^{_{\rm sc}}_{_{P,Q}}(1)$ and denoting by  $ I_Q$ the {\it finite\/} set of pairs $(Q', a')$ formed by an $\F\-$selfcentralizing subgroup $Q'$ of $P$ and by a $\P^{^{\rm sc}}\!\-$morphism 
$a'\,\colon Q'\to Q$ belonging to~$\P^{^{\rm sc}}\! (Q,Q')_{i_{Q'}^P}$ [12,~4.5.1], we may assume that
$$Q\,\hat\times\, P = \bigoplus_{(Q',a')\in  I_Q} Q'
\eqno £2.15.10.$$
Actually, it follows from [12,~4.10] that we have a functor 
$ I\,\colon \P^{^{\rm sc}}\!\to \Set$  mapping any 
$\F\-$selfcentralizing subgroup $Q$ of $P$ on the {\it finite\/} set  $ I_Q$ 
and any $\P^{^{\rm sc}}\!\-$morphism $ x\,\colon R\to Q$ on the map  
$ I_{x}\,\colon I_R\to I_Q$ determined by the $\ad (\P^{^{\rm sc}})\-$morphism
$$x\,\hat\times\,  i_P^P : R\,\hat\times\, P\too Q\,\hat\times\, P
\eqno £2.15.11;$$
that is to say, any element $(R',b')$ in $ I_R$ determines an element  
$(Q',a')$ in~$ I_Q$ in such a way that $Q'$ contains $R'$ and that, setting 
$ i_{R'}^{Q'} = \tau^{_{\rm sc}}_{_{Q',R'}}(1)\,,$ we have
$$a'\,\.\, i_{R'}^{Q'}  = x\, \.\, b'
\eqno £2.15.12.$$

\smallskip
Now, in order to  prove that $\Bbb H_{\I }^n \big(\P^{^{\rm sc}}\!, 
\frak h_0(\bar\frak m^{_\O})\big) = \{0\}$ for any $n\ge 1\,,$ we will quote [12,~Theorem~3.5]; for this purpose,
we need to consider a {\it homotopic system\/}~$\H$ --- as introduced in [12,~2.6] --- associated with $\T^{^{\rm sc}}_{\!P}\,,$ with the 
{\it  $\T^{^{\rm sc}}_{\!P}\-$category\/}~$\P^{^{\rm sc}}\!$ and with the  subcategory $\I$ of $\P^{^{\rm sc}}\!$; our {\it homotopic system\/} $\H$ is the quintuple formed by the {\it interior structure  $\I$\/} above, by the {\it trivial co-interior structure\/} of $\P^{^{\rm sc}}\!$, by the functor 
$ I\,\colon \P^{^{\rm sc}}\!\to \Set$ above, by the functor [12,~2.5]
$$\frak w : I\rtimes \P^{^{\rm sc}}\!\too 
\widetilde\F^{^{\rm sc}}_{\!P}\i \widetilde\F^{^{\rm sc}}\!
\eqno £2.15.13\phantom{.}$$
mapping any $ I\rtimes \P^{^{\rm sc}}\!\-$object $(Q',a',Q)$ on $Q'$ 
and any $ I\rtimes \P^{^{\rm sc}}\!\-$morphism 
$$(i_{R'}^{Q'},x) : (R',b',R)\too (Q',a',Q)
\eqno £2.15.14\phantom{.}$$
on $\tilde\iota_{R'}^{Q'} \,\colon R'\to Q'\,,$ where $Q$ and $R$ are $\F\-$selfcentralizing subgroups of $P\,,$  $(Q',a')$ and $(R', b')$ are respective elements of $ I_Q$ and $ I_R\,,$ and $x\,\colon R\to Q$ is a $\P^{^{\rm sc}}\!\-$morphism fulfilling equality~£2.15.12 above; finally, denoting by 
$$\tilde\frak p :   I\rtimes \P^{^{\rm sc}}\! \too \widetilde\F^{^{\rm sc}}\!
\eqno £2.15.15\phantom{.}$$
the {\it forgetful\/} functor mapping $(Q',a',Q)$ on $Q$ and $(i_{R'}^{Q'},x))$ on the image $\tilde x$ of~$x$ in 
$\widetilde\F^{^{\rm sc}}\!(Q,R)\,,$ the fifth term
in $\H$ is the {\it natural map\/}
$$\omega : \frak w\too \tilde\frak p
\eqno £2.15.16\phantom{.}$$
sending any $ I\rtimes \P^{^{\rm sc}}\!\-$object $(Q',a',Q)$ to the image
$$\matrix{\tilde a' : \hskip-20pt&Q' \hskip-10pt&\too &\hskip-10ptQ\cr 
&\Vert\hskip-10pt&\phantom{\Big\uparrow}&\hskip-10pt\Vert\cr
& \frak w (Q',a',Q)\hskip-10pt&
&\hskip-10pt\tilde\frak p (Q',a',Q)\cr}
\eqno £2.15.17\phantom{.}$$ 
of $a'$ in $\widetilde\F^{^{\rm sc}}\! (Q,Q')\,;$ the {\it naturalness\/} of $\omega$ is then easily checked from equality~£2.15.12.

\smallskip
At this point, following [12,~2.9 and~2.10], from the {\it homotopic system\/}~$\H$ and from the  {\it contravariant\/} 
functor $\frak h_0(\bar\frak m^{_{\O}})\,,$ which factorizes through the canonical 
functor~$\P^{^{\rm sc}}\!\to \widetilde\F^{^{\rm sc}}\!$, we get a {\it contravariant\/} functor
$$\H \big(\frak h_0(\bar\frak m^{_{\O}})\big) : \P^{^{\rm sc}}\too  k\-\mod
\eqno £2.15.18\phantom{.}$$
 sending any $\F\-$selfcentralizing subgroup $Q$ of $P$ to
 $$\Big(\H \big(\frak h_0(\bar\frak m^{_{\O}})\big)\Big)(Q) = 
 \big(\prod_{(Q',a')\in  I_Q} (\bar M^{^{\O}})^{Q'}\big)^Q
 \eqno £2.15.19,$$
 and a {\it natural map\/}
 $$\Delta_\H \big(\frak h_0(\bar\frak m^{_{\O}})\big) : \frak h_0 
 (\bar\frak m^{_{\O}})\too \H \big(\frak h_0(\bar\frak m^{_{\O}})\big) 
 \eqno £2.15.20\phantom{.}$$
 sending any $\F\-$selfcentralizing subgroup $Q$ of $P$ to the $ k\-$module
 homomorphism
 $$\Delta_\H \big(\frak h_0(\bar\frak m^{_{\O}})\big)_Q : (\bar M^{^{\O}})^{Q} \too \big(\prod_{(Q',a')\in \hat I_Q} (\bar M^{^{\O}})^{Q'}\big)^Q
\eqno £2.15.21\phantom{.}$$ 
mapping any $\bar m\in  (\bar M^{^{\O}})^{Q}$ on 
$$\Delta_\H \big(\frak h_0(\bar\frak m^{_{\O}})\big)_Q (\bar m) = 
\sum_{(Q',a')\in  I_Q} \big(\bar\frak m^{_{\O}} (a')\big)(\bar m)
\eqno £2.15.22\,.$$

\smallskip
Then, it follows from [12,~Theorem~3.5] that, for our purpose, it suffices to exhibit a {\it natural section\/} of $\Delta_\H \big(\frak h_0
(\bar\frak m^{_{\O}})\big)$
$$\theta : \H \big(\frak h_0(\bar\frak m^{_{\O}})\big)\too
 \frak h_0(\bar\frak m^{_{\O}}) 
 \eqno £2.15.23;$$
  \eject
 \noindent
 explicitly, for any $\F\-$selfcentralizing subgroup $Q$ of $P$ we will define a section 
$$\theta_Q : \big( \prod_{(Q',a')\in  I_Q} (\bar M^{^{\O}})^{Q'}\big)^Q\too 
(\bar M^{^{\O}})^Q
\eqno £2.15.24\phantom{.}$$ 
of $\Delta_\H \big(\frak h_0(\bar\frak m^{_{\O}})\big)_Q\,.$ Note that
the action of $u\in Q$ maps the element
$$\bar m = \sum_{(Q',a')\in I_Q} \bar m_{(Q',a')}\in \prod_{(Q',a') \in  
I_Q} (\bar M^{^{\O}})^{Q'}
 \eqno £2.15.25\phantom{.}$$
on $\sum_{(Q',a')\in I_Q} \bar m_{(Q',\tau^{_{\rm sc}} (u)\.a')}\,;$ thus, 
$\bar m$ belongs to $ \big({\displaystyle \prod_{(Q',a')\in  I_Q}}
(\bar M^{^{\hat\O}})^{Q'}\big)^Q$ if and only if we have 
$\bar m_{(Q',\tau^{_{\rm sc}} (u)\.\skew2\hat{\hat a}')} = \bar m_{(Q',a')}$ for any $u\in Q\,;$ that is to say, $\bar m_{(Q',a')}$ only depends on the 
pair~$(Q',\tilde a')\,.$

\smallskip
Actually, denoting by $\widetilde\I_Q$ the set of pairs $(Q',\tilde a')$ when 
$(Q',a')$ runs over~$\I_Q\,,$ it follows from [12,~Corollary~4.7] that the 
{\it direct product\/} in $\ad (\P^{^{\rm sc}})$ induces a {\it direct product\/}
 in $\ad (\widetilde\F^{^{\rm sc}})$ --- noted  $\hat{\widetilde\times}$ --- and that we may assume that
$$Q\,\hat{\widetilde\times}\, P = \bigoplus_{(Q',\tilde a')\in \widetilde\I_Q} Q'
\eqno £2.15.26;$$
we denote by $\widetilde\I_Q^\circ$ the set of pairs $(Q',\tilde a')\in \widetilde\I_Q$ --- called  {\it extremal\/} --- where $\tilde a'$ is an isomorphism and it is easily checked that we have a canonical bijection
$$\widetilde\I_Q^\circ\cong \widetilde\F (P,Q)
\eqno £2.15.27.$$
With all this notation, for any element $\bar m = 
\sum_{(Q',a')\in  I_Q} \bar m_{(Q',a')}$ belonging to
$$\Big(\H \big(\frak h_0(\bar\frak m^{_{\O}})\big)\Big)(Q) = 
\big( \prod_{(Q',a')\in  I_Q} (\bar M^{^{\O}})^{Q'}\big)^Q
\eqno £2.15.28,$$
we define
$$\theta_Q (\bar m ) ={1\over \big\vert \widetilde\F^{^{\rm sc}} (P,Q)\big\vert}\sum_{(Q',\tilde a')\in \widetilde\I^\circ_Q} 
\big(\bar\frak m^{_{\O}} (a')\big)^{-1} 
(\bar m_{(Q',a')})
\eqno £2.15.29\phantom{.}$$
where, for any $(Q',\tilde a')\in \widetilde\I^\circ_Q\,,$ $a'$ denotes a representative of $\tilde a'$ in 
$\P^{^{\rm sc}}(Q,Q')\,;$ this makes sense since it follows from [9,~6.6.4 and Proposition~6.7] and from condition~£1.3.2 above that $p$ does not divide $\vert \widetilde\F (P,Q)\vert$ and, since the element 
$\big(\bar\frak m^{_{\O}} (a')\big)^{-1} (\bar m_{(Q',a')})$ belongs to 
$(\bar M^{^{\O}})^{Q}\,,$ it is clear that $\theta_Q (\bar m )$ does not depend on the choice of the representative $a'\,;$ moreover, it follows from definition~£2.15.22 above that $\theta_Q$ is indeed a section of~$\Delta_\H \big(\frak h_0 (\bar\frak m^{_{\O}})\big)_Q\,.$

\smallskip
It only remains to prove that the correspondence sending any $\F\-$selfcen-tralizing subgroup $Q$ of $P$ to $\theta_Q$ is {\it natural\/}; that is to say,
for any $\P^{^{\rm sc}}\!\-$mor-phism $x\,\colon R\to Q$ we have 
 to prove the commutativity of the
following diagram
$$\matrix{\Big(\H \big(\frak h_0(\bar\frak m^{_{\O}})\big)\Big)(Q)
&\buildrel \theta_Q\over \too
& \bar\frak m^{_{\O}} (Q)\cr
\hskip-70pt{\scriptstyle (\H (\frak h_0(\bar\frak m^{_{\O}})))(x)}\hskip5pt \big\downarrow
&\phantom{\Big\downarrow}&\big\downarrow\hskip5pt {\scriptstyle
\bar\frak m^{_{\O}} (x)}\hskip-20pt\cr
\Big(\H \big(\frak h_0(\bar\frak m^{_{\O}})\big)\Big)(R)
&\buildrel \theta_R\over \too &\bar\frak m^{_{\O}} (R)\cr}
\eqno £2.15.30.$$
Explicitly, it follows from~£2.15.11 and~£2.15.12 above that the $ k\-$module
homomorphism
$$\Big(\H \big(\frak h_0(\bar\frak m^{_{\O}})\big)\Big)(x) : 
\Big(\H \big(\frak h_0(\bar\frak m^{_{\O}})\big)\Big)(Q)\too
\Big(\H \big(\frak h_0(\bar\frak m^{_{\O}})\big)\Big)(R)
\eqno £2.15.31\phantom{.}$$
sends the element ${\bar m }= \sum_{(Q',a')\in  I_Q} \bar m_{(Q',a')}$ above to 
$$\sum_{(R',b')\in  I_R} \big(\bar\frak m^{_{\O}}( i^{Q'}_{R'})\big) 
({\bar m}_{(Q',a')})
\eqno £2.15.32\phantom{.}$$
where, for any $(R',b')\in  I_R\,,$ $(Q',a')$ is the unique element of $ I_Q$ such that $R'\i Q'$ and $a'\,\.\,  i_{R'}^{Q'}  = x\, \.\,b'$ (cf.~£2.15.12); then, $\theta_R$ maps this element on
$${1\over \vert \widetilde\F^{^{\rm sc}} (P,R)\vert} \sum_{(R',\tilde b')\in 
\widetilde\I^\circ_R} \big(\bar\frak m^{_{\O}} (b')\big)^{-1} 
\Big(\big(\bar\frak m^{_{\O}}( i^{Q'}_{R'})\big) (\bar m_{(Q',a')})\Big)
\eqno £2.15.33\phantom{.}$$
where, for any $(R',\tilde b')\in \widetilde\I^\circ_R\,,$ $b'$ denotes a representative of 
$\tilde b'$ in $\P^{^{\rm sc}}\!(R,R')\,.$

\smallskip
But, it follows from [12,~Lemma~4.4] that the category $\ad (\widetilde\F^{^{\rm sc}})$ admits {\it pull-backs\/} and, more precisely, that we have the {\it pull-back\/}
$$\matrix{&&Q\cr
&{\tilde x\atop}\hskip-5pt\nearrow&&\hskip-30pt\nwarrow\cr
R&&&&\hskip-25pt Q\,\hat{\widetilde\times}\, P\cr
&\nwarrow&&\nearrow \hskip-5pt{\atop \tilde x\,\hat{\widetilde\times}\, i_P^P }\cr
&&R\,\hat{\widetilde\times}\, P}
\eqno £2.15.34;$$
explicitly, for any $(Q',\tilde a')\in \widetilde\I_Q\,,$ choosing a representative $a'$ of $\tilde a'$ and  a set of representatives 
$W_{\!(Q',\tilde a')}$ for the set of double classes $\varphi_{x} (R)\backslash Q /\varphi_{a'}(Q')$ where 
$\varphi_x\in \F (Q,R)$ and $\varphi_{a'}\in \F (Q,Q')$ are the respective images of $x$ and~$a'\,,$ it is well-known that we get  the {\it pull-back\/} in  $\ad (\widetilde\F^{^{\rm sc}})$
$$\matrix{&&Q\cr
&{\tilde x\atop}\hskip-5pt\nearrow\hskip-20pt&
&\hskip-20pt\nwarrow \hskip-5pt{\tilde a'\atop}\cr
R\hskip-10pt&&&&\hskip-10pt Q'\cr
&\nwarrow\hskip-20pt&&\hskip-20pt\nearrow \cr
&&{\displaystyle\bigoplus_{w\in W_{\!(Q',\tilde a')}^{^{\rm sc}}}} R'_w\cr} 
\eqno £2.15.35\phantom{.}$$
where $W_{\!(Q',\tilde a')}^{^{\rm sc}}$ is the subset of 
$w\in W_{\!(Q',\tilde a')}$ such that the subgroup of $Q'$
$$R'_w =  \varphi_{a'}^{-1} \big(\varphi_x (R)^w\big)
\eqno £2.15.36\phantom{.}$$
is still $\F\-$selfcentralizing; moreover, denote by $b'_w$ the element in $\P^{^{\rm sc}}(R,R'_w)$ fulfilling 
$a'\,\.\, i_{R'_w}^{Q'}  = x\, \.\,b'_w\,.$ At this point, in the category 
 $\ad (\widetilde\F^{^{\rm sc}})$ we get
$$R\,\hat{\widetilde\times}\, P\cong \bigoplus_{(Q',\tilde a')\in \widetilde\I_Q}\, 
\bigoplus_{w\in W_{\!(Q',a')}^{^{\rm sc}}} R'_w
\eqno £2.15.37\phantom{.}$$
and we actually may assume that
$$\widetilde\I_R = \bigsqcup_{(Q',\tilde a')\in \widetilde\I_Q}\{(R'_w,\tilde b'_w)\}_{w\in 
W_{\!(Q',\tilde a')}^{^{\rm sc}}}
\eqno £2.15.38.$$

\smallskip
Consequently, denoting by $W_{\!(Q',\tilde a')}^\circ$  the subset of $w\in W_{\!(Q',\tilde a')}$ such that $R'_w$ is 
$\F\-$isomorphic to $R\,,$ we clearly have $W_{\!(Q',\tilde a')}^\circ\i W_{\!(Q',\tilde a')}^{^{\rm sc}}$ and
from~£2.15.32 and~£2.15.33 we get
$$\eqalign{&\big\vert \widetilde\F^{^{\rm sc}} (P,R)\big\vert\.\bigg(\theta_R \circ 
\Big(\H \big(\frak h_0(\bar\frak m^{_{\O}})\big)\Big)(x)\bigg) (\bar m)\cr
& = \sum_{(Q',\tilde a')\in \widetilde\I_Q}\, \sum_{w\in W_{\!(Q',\tilde a')}^\circ}
\big(\bar\frak m^{_{\O}} (b'_w)^{-1} 
\circ \bar\frak m^{_{\O}}( i^{Q'}_{R'_w})\big) (\bar m_{(Q',a')})\cr}
\eqno £2.15.39\phantom{.}$$
where, for any $(Q',\tilde a')\in \widetilde\I_Q$ and any  $w\in W_{\!(Q',\tilde a')}^\circ\,,$ we have $a'\,\.\,  i_{R'_w}^{Q'}  = 
x\, \.\, b'_w\,;$ in particular, we still have
 $$\bar\frak m^{_{\O}}( i^{Q'}_{R'_w}) \circ \bar\frak m^{_{\O}} ( a') =
\bar\frak m^{_{\O}} (b'_w) \circ \bar\frak m^{_{\O}} ( x)
 \eqno £2.15.40\phantom{.}$$
 and therefore the composition
 $$\bar\frak m^{_{\O}} (b'_w)^{-1} 
\circ \bar\frak m^{_{\O}}(i^{Q'}_{R'_w}) = 
\bar\frak m^{_{\O}} (x) \circ
\bar\frak m^{_{\O}} ( a')^{-1}
\eqno £2.15.41\phantom{.}$$
does not depend on $w\in W_{\!(Q',\tilde a')}^\circ\,.$

\smallskip
On the other hand, from definition~£2.15.29 we get
$$\eqalign{\big\vert \widetilde\F^{^{\rm sc}} (P,Q)\big\vert\. &\big(\bar\frak m^{_{\O}} (x) \circ \theta_Q\big) (\bar m)\cr
& = \sum_{(Q',\tilde a')\in \widetilde\I^\circ_Q} \big(\bar\frak m^{_{\O}} (x) \circ \bar\frak m^{_{\O}} (a')^{-1} \big)
(\bar m_{(Q',a')})\cr}
\eqno £2.15.42.$$
Hence, in order to prove the commutativity of diagram~£2.15.30, it suffices to show that, for any 
$(Q',\tilde a')\in \widetilde\I_Q\,,$ either $(Q',\tilde a')$ belongs to  $\widetilde\I^\circ_Q$ and we have $\vert W_{\!(Q',\tilde a')}^\circ\vert = 1\,,$ or otherwise  $p$ divides $\vert W_{\!(Q',\tilde a')}^\circ\vert\,;$ but, the set
$W_{\!(Q',\tilde a')}^\circ$ is actually a set of representatives for 
$\varphi_{a'} (Q')\backslash \T_Q\big(\varphi_{a'} (Q'),\varphi_x (R)\big)$ and if this quotient is not empty then $N_{Q}\big(\varphi_{a'}(Q')\big)/\varphi_{a'} (Q')$ acts freely on it. We are done.
\eject

\bigskip
\noindent
{\bf Corollary~£2.16.} {\it  For any pair of $\K_*\hat\P^{^{\rm sc}}\!\-$modules 
$\frak m$ and~$\frak m'$ we have
$$\langle X_\frak m,X_{\frak m'}\rangle = 
{\rm dim}_\K \big(\Nat (\frak m',\frak m)\big)
\eqno £2.16.1.$$\/}

\par
\noindent
{\bf Proof:} According to~£2.8.4 we have 
$$\Nat (\frak m',\frak m) = {\displaystyle \lim_{\longleftarrow}}\,
(\frak m^*\otimes_\K \frak m') = \Bbb H_{\rm r}^0  (\P^{^{\rm sc}}\!,
\frak m^*\otimes_\K \frak m')
\eqno £2.16.2\phantom{.}$$
 and therefore,  setting $C^n = \Bbb C_{\rm r}^n  (\P^{^{\rm sc}}\!,\frak m^*\otimes_\K \frak m')$
for any $n\in \Bbb N$ (cf.~£2.6), it follows from Theorem~£2.15 that  we have an infinite exact sequence
$$0\too \Nat (\frak m',\frak m) \too C^0 \too \dots\too C^n\too C^{n+1}\too \dots
\eqno £2.16.3;$$
actually,  we can identify $C^n$ with the set of elements
$$(m_{\frak q})_{\frak q\in \Fct_{\rm r} (\Delta_n,\P^{^{\rm sc}})} \in 
\prod_{\frak q\in \Fct_{\rm r} (\Delta_n,\P^{^{\rm sc}}\!)} (\frak m^*\otimes_\K \frak m')\big(\frak q (0)\big)
\eqno £2.16.4\phantom{.}$$
such that, for any natural isomorphism $ x\,\colon \frak q\cong \frak q'$ between {\it regular $\P^{^{\rm sc}}\!\-$valued 
$n\-$chains\/} $\frak q$ and $\frak q'\,,$ $(\frak m^*\otimes_\K \frak m') ( x_0)$ maps  $m_{\frak q'}$ on $m_{\frak q}\,;$
that is to say, we actually have
$$\Bbb C_{\rm r}^n  (\P^{^{\rm sc}}\!,\frak m^*\otimes_\K \frak m')
\cong \prod_{\frak q} (\frak m^*\otimes_\K \frak m')\big(\frak q (0)\big)^{\P^{^{\rm sc}}\!(\frak q)}
\eqno £2.16.5\phantom{.}$$
where $\frak q$ runs over a set of representatives for the set of isomorphism
classes in $\Fct_{\rm r} (\Delta_n,\P^{^{\rm sc}})$
[9,~£A5.3].

\smallskip
On the other hand, it is clear that for $n$ big enough there are no {\it regular $\P^{^{\rm sc}}\!\-$valued $n\-$chains\/} and therefore, in the exact sequence above, only finitely many terms are not zero; thus, we still get
$${\rm dim}_K \big(\Nat (\frak m',\frak m)\big) = \sum_{(\frak q,\Delta_n)} (-1)^{n}\,{\rm dim\,}_\K 
\Big((\frak m^*\otimes_\K \frak m') \big(\frak q (0)\big)^{\P^{^{\rm sc}}\!(\frak q)}\Big)
\eqno £2.16.6\phantom{.}$$
where $(\frak q,\Delta_n)$ runs over a set of representatives for the isomorphism classes of {\it regular\/} 
$\ch^*(\P^{^{\rm sc}} )\-$objects (cf.~£A5.3). Moreover, for any functor $\hat\frak q\,\colon\Delta_n\to \hat\P^{^{\rm sc}} $ lifting $\frak q$ we have
$$(\frak m^*\otimes_\K \frak m')\big(\frak q (0)\big) = {\rm Hom}_\K \Big(\frak m \big(\hat\frak q (0)\big),
\frak m' \big(\hat\frak q (0)\big)\Big)
\eqno £2.16.7\phantom{.}$$
and, in particular, we get
$$(\frak m^*\otimes_\K \frak m')\big(\frak q (0)\big)^{\P^{^{\rm sc}}\!(\frak q)} =  
{\rm Hom}_{\K _*\hat\P^{^{\rm sc}}\!(\hat\frak q)}
\Big(\frak m \big(\hat\frak q (0)\big),\frak m' \big(\hat\frak q (0)\big)\Big)
\eqno £2.16.8;$$
thus, denoting by $X_{\frak m (\hat\frak q (0))}$ and $X_{\frak m' (\hat\frak q (0))}$ the respective classes of 
$\frak m \big(\hat\frak q (0)\big)$ and $\frak m' \big(\hat\frak q (0)\big)$ in 
$\G_\K \big(\hat\P^{^{\rm sc}}\!(\hat\frak q)\big)\,,$ we obtain 
(cf.~£2.1 and~£2.4)
$${\rm dim\,}_\K \Big((\frak m^*\otimes_\K \frak m')
\big(\frak q (0)\big)^{\P^{^{\rm sc}}\!(\frak q)}\Big) = 
\langle X_{\frak m (\hat\frak q (0))},X_{\frak m' (\hat\frak q (0))} \rangle
\eqno £2.16.9.$$
Now, equality~£2.16.1 follows from equality~£2.16.6.
\eject

\bigskip
\noindent
{\bf Corollary~£2.17.} {\it The $\O\-$module homomorphism 
$$\cat_\K : \G (\K_*\hat\P^{^{\rm sc}}\!\-\mod)\too
 \G_\K (\hat\P^{^{\rm sc}})
 \eqno £2.17.1\phantom{.}$$
 is injective. In particular, there are finitely many isomorphism classes of simple
 $\K_*\hat\P^{^{\rm sc}}\!\-$modules.\/}
 
  \medskip
 \noindent
 {\bf  Proof:} It is clear that if $\frak m$ and $\frak m'$ are nonisomorphic  simple 
 $\phantom{.}\hat\P^{^{\rm sc}}\!\-$modules then we have 
 $\Nat (\frak m,\frak m') =\{0\}\,;$ consequently, if $\{X_i\}_{i\in I}$ is a finite family in $\phantom{.}\G (\K_*\hat\P^{^{\rm sc}}\!\-\mod)$
 of classes of   simple $\K_*\hat\P^{^{\rm sc}}\!\-$modules and for a family 
 $\{\lambda_i\}_{i\in I}$ in~$\O$ we have $\sum_{i\in I} \lambda\.\cat_\K (X_i) = 0\,,$ it suffices to perform the scalar product by $\cat_\K (X_j)$ to get $\lambda_j = 0$ for any 
 $j\in I\,.$

 \medskip
 £2.18. In order to prove that $\cat_\K\,\colon \G (\K_*\hat\P^{^{\rm sc}}\!\-\mod)\to \G_\K (\hat\P^{^{\rm sc}})$ is also surjective, following Lemma~£1.6 above we consider a suitable finite subfield $\hat k$ of~$k$ and a regular central 
$\hat k^*\-$extension $\skew4\hat{\hat\F}^{^{\rm sc}}$ of~$\F^{^{\rm sc}}$ such that the extension of 
$\skew4\hat{\hat\F}^{^{\rm sc}}$ from $\hat k^*$ to $k^*$ is equivalent to $\hat\F^{^{\rm sc}}\!\,;$ denote by $\hat\O$ the converse image of $\hat k$ in $\O$ and by~$\hat\K$ the {\it field of quotients\/} of $\hat\O\,.$ Then, denoting by $\skew4\hat{\hat\P}^{^{\rm sc}}\!$ the converse image of~$\skew4\hat{\hat\F}^{^{\rm sc}}$ in $\hat\P^{^{\rm sc}}$,
we may choose $\hat k$  big enough to get
$$\G_{\hat\K} \big(\skew4\hat{\hat\P}^{^{\rm sc}} (\hat\frak q)\big) \cong 
\G_\K \big(\hat\P^{^{\rm sc}} (\hat\frak q)\big)
\eqno £2.18.1\phantom{.}$$
for~any  $\skew4\hat{\hat\P}^{^{\rm sc}}\!\-$ chain $\hat\frak q\,.$

\medskip
£2.19.  Let $\ell$ be a prime number not dividing neither $\vert\hat k^*\vert$  nor $\vert \F (Q) \vert$ for 
any~$\F\-$selfcentralizing subgroup $Q$ of $P\,,$ and denote by $\O_\ell$ a complete discrete~valuation ring with a
{\it quotient field\/} $\K_\ell$ of characteristic zero and a {\it finite residue field\/} $k_\ell$ of characteristic $\ell\,;$ 
we can choose $\O_\ell$   in such a way that $(k_\ell)^*$ would contain a (unique) subgroup isomorphic to $\hat k^*\,;$ 
in particular,  choosing~an inclusion $\hat k^*\i (k_\ell)^*\,,$ any $\hat k^*\-$group $\skew4\hat{\hat G}$ induces a 
$(k_\ell)^*\-$group  $(k_\ell)^*\times_{\hat k^*} \skew4\hat{\hat G}$ and, since this correspondence is functorial, 
from $\skew4\hat{\hat\P}^{^{\rm sc}}\!$ we actually get a  $(k_\ell)^*\-$category $\skew4\hat{\hat\P}^{^{\rm sc,\ell}}$
containing $\skew4\hat{\hat\P}^{^{\rm sc}}\,.$ Moreover, {\it via\/} a suitable field $\skew2\hat{\hat\K}$ containing $\hat\K$ and $\K_\ell\,,$ it is clear that choosing $\O_\ell$ big enough for~any  {\it $\skew4\hat{\hat\P}^{^{\rm sc}}\!\-$chain\/} 
$\hat\frak q\,\colon \Delta_n\to \skew4\hat{\hat\P}^{^{\rm sc}}\!$ we can get an $\O\-$module isomorphism
$$\G_{\hat\K} \big(\skew4\hat{\hat\P}^{^{\rm sc}} (\hat\frak q)\big) 
\cong \G_{\K_\ell} \big(\skew4\hat{\hat\P}^{^{\rm sc,\ell}} (\hat\frak q^{^\ell})\big)
\eqno £2.19.1\phantom{.}$$
where $\hat\frak q^{^\ell}$ is determined by $\hat\frak q$ and by the inclusion 
$\skew4\hat{\hat\P}^{^{\rm sc}}\i \skew4\hat{\hat\P}^{^{\rm sc,\ell}}\,.$ Finally, according to our choice of~$\ell\,,$ we know that the {\it Brauer decomposition map\/} determines  an $\O\-$module isomorphism
$$\G_{\K_\ell} \big(\skew4\hat{\hat\P}^{^{\rm sc,\ell}} (\hat\frak q^{^\ell})\big)
\cong \G_{k_\ell} \big(\skew4\hat{\hat\P}^{^{\rm sc,\ell}} (\hat\frak q^{^\ell})\big)
\eqno £2.19.2$$
 for~any  {\it $\skew4\hat{\hat\P}^{^{\rm sc}}\!\-$chain\/} $\hat\frak q\,\colon \Delta_n\to 
 \skew4\hat{\hat\P}^{^{\rm sc}}$.
 \eject

\medskip
£2.20. Consequently,  for~any  {\it $\skew4\hat{\hat\P}^{^{\rm sc}}\!\-$chain\/} $\hat\frak q\,\colon \Delta_n\to 
 \skew4\hat{\hat\P}^{^{\rm sc}}$, the composition of all these $\O\-$module isomorphisms
supplies an $\O\-$module isomorphism
$$\G_\K \big(\hat\P^{^{\rm sc}} (\hat\frak q)\big)\cong 
\G_{k_\ell} \big(\skew4\hat{\hat\P}^{^{\rm sc,\ell}} (\hat\frak q^{^\ell})\big)
\eqno £2.20.1\phantom{.}$$
and therefore, since all of them are functorial, we get a {\it natural isomorphism\/}
$$\frak g_\K\circ \aut_{\hat\P^{^{\rm sc}}}\cong \frak g_{k_\ell}
\circ \aut_{\skew4\hat{\hat\P}^{^{\rm sc,\ell}} }
\eqno £2.20.2,$$
so that we still get an $\O\-$module isomorphism
$$\G_\K (\hat\P^{^{\rm sc}}) = 
\lim_{\longleftarrow}\,(\frak g_\K\circ \aut_{\hat\P^{^{\rm sc}}}) \cong
\lim_{\longleftarrow}\,(\frak g_{k_\ell}\circ 
\aut_{\skew4\hat{\hat\P}^{^{\rm sc,\ell}} }) = 
\G_{k_\ell} (\skew4\hat{\hat\P}^{^{\rm sc,\ell}} )
\eqno £2.20.3.$$
Moreover note that, as in~£1.9 above, the {\it faithful\/} functor $\tau^{_{\rm sc}}\,\colon 
\T^{^{\rm sc}}_P \to \P^{^{\rm sc}}$ can be lifted to a {\it faithful\/} functor 
$\skew2\hat{\hat\tau}^{^{\rm sc}}\,\colon \T^{^{\rm sc}}_P \to \skew4\hat{\hat\P}^{^{\rm sc}} \i
\skew4\hat{\hat\P}^{^{\rm sc,\ell}}$ [13,~Proposition~3.5], and we set 
$\hat i^Q_R =\skew2\hat{\hat\tau}^{_{\rm sc}}_{Q,R}(1)$ for any pair of $\F\-$selfcentralizing subgroups $Q$ and $R$ of $P$ fulfilling $R\i Q\,;$ in particular, $\skew4\hat{\hat\P}^{^{\rm sc}}$ and $\skew4\hat{\hat\P}^{^{\rm sc,\ell}}$ become {\it divisible $\F^{^{\rm sc}}\-$localities\/} and therefore, for another such a pair $Q'$ and $R'\,,$ we have a {\it restriction\/} map
$$\hat r_{R',R}^{Q',Q} : \skew4\hat{\hat\P}^{^{\rm sc,\ell}}\! (Q',Q)_{R',R}
\too \skew4\hat{\hat\P}^{^{\rm sc,\ell}} (R',R)
\eqno £2.20.4$$
fulfilling $\hat i^{Q'}_{R'}\,\.\,\hat r_{R',R}^{Q',Q} (\hat x) = \hat x\,\.\, \hat i^Q_R$
for any $\hat x\in \skew4\hat{\hat\P}^{^{\rm sc,\ell}}\! (Q',Q)$ sending $R$ to $R'\,.$

£2.21. At this point, let us  call {\it reduced $(k_\ell)_*\skew4\hat{\hat\P}^{^{\rm sc,\ell}}\!\-$module\/} any  {\it contra-variant $(k_\ell)^*\-$functor\/}~$\frak m_\ell\,\colon \skew4\hat{\hat\P}^{^{\rm sc,\ell}}\!\to k_\ell\-\mod$  mapping any 
$\F\-$selfcentralizing subgroup $Q$ of~$P$ on the {\it same\/} finite dimensional $k_\ell\-$module $M_\ell$ and, for any 
$\F\-$selfcentralizing subgroup $R$ of~$Q\,,$ the $\skew4\hat{\hat\P}^{^{\rm sc,\ell}}\!\-$morphism $\hat i^Q_R$ on ${\rm id}_{M_\ell}\,;$ this is coherent with our terminology above since the group ${\rm GL}_{k_\ell} (M_\ell)$ and therefore the 
$(k_\ell)^*\-$subgroup $\hat G(\frak m_\ell)$ of 
${\rm GL}_{k_\ell} (M_\ell)$ generated  by~$\frak m_\ell (\hat x)\,,$ where  
$\hat x\,\colon R\to Q$ runs over the set of  
$\skew4\hat{\hat\P}^{^{\rm sc,\ell}}\!\-$morphisms, are {\it finite\/}.

 \bigskip
 \noindent
 {\bf Theorem~£2.22.} {\it With the notation above, let $\{X_{\hat\frak q}\}_{\hat\frak q}$ where  $\hat\frak q$ 
 runs over the set of  $\,\skew4\hat{\hat\P}^{^{\rm sc,\ell}}\!\!\-$chains be a family which  belongs to 
 $\G_{k_\ell} (\skew4\hat{\hat\P}^{^{\rm sc,\ell}})$ and fulfills that, for such a 
 $\hat\frak q\,,$ $X_{\hat\frak q}$ is the class of a $(k_\ell)_*\skew4\hat{\hat\P}^{^{\rm sc,\ell}}\!\-$module 
 $M_{\hat\frak q}\,.$ Then, there exists a reduced $(k_\ell)_*\skew4\hat{\hat\P}^{^{\rm sc,\ell}}\!\-$module
$\frak m_\ell\,\colon \skew4\hat{\hat\P}^{^{\rm sc,\ell}}\!\to k_\ell\-\mod$  
such that $X_{\hat\frak q}$ is the class of $\frak m_\ell \big(\hat\frak q (0)\big)$
in $\G_{k_\ell} \big(\skew4\hat{\hat\P}^{^{\rm sc,\ell}} (\hat\frak q)\big)$
for~any  $\skew4\hat{\hat\P}^{^{\rm sc,\ell}}\! \!\-$chain $\hat\frak q\,.$\/}
 
 \medskip
 \noindent
 {\bf  Proof:} Let $\frak X$ be a nonempty set  of  $\F\-$selfcentralizing  subgroups of $P$ which contains any subgroup of $P$ admitting an $\F\-$morphism from some subgroup in~$\frak X\,,$ and respectively denote by 
$\T^{^{\frak X}}_P\,,$ $\F^{^\frak X}\,,$  $\P^{^\frak X}$ and $\skew4\hat{\hat\P}^{^{\frak X,\ell}}$ the {\it full\/} subcategories of  $\T^{^{\rm sc}}_P\,,$  $\F^{^{\rm sc}}\,,$  $\P^{^{\rm sc}}$ and 
 $\skew4\hat{\hat\P}^{^{\rm sc,\ell}}$ over $\frak X$ as the  set of objects; arguing by induction on~$\vert\frak X\vert\,,$ we  prove that there is a  {\it reduced $(k_\ell)_*\skew4\hat{\hat\P}^{^{\frak X,\ell}} \!\!\-$module\/}  
 $\frak m_\ell^{_{\frak X}}\,\colon \skew4\hat{\hat\P}^{^{\frak X,\ell}} \!\to  
 k_\ell\-\mod$ such that $X_{\hat\frak q}$ is the class of 
 $\frak m_\ell \big(\hat\frak q (0)\big)$
in $\G_{k_\ell} \big(\skew4\hat{\hat\P}^{^{\frak X,\ell}} (\hat\frak q)\big)$
for~any  {\it $\skew4\hat{\hat\P}^{^{\frak X,\ell}} \!\-$chain\/} $\hat\frak q\,.$
\eject

\smallskip
If $\frak X = \{P\}$ then $\skew4\hat{\hat\P}^{^{\frak X,\ell}} \!$ has just one object $P\,;$ in this situation, still denoting by 
$P$ the $\skew4\hat{\hat\P}^{^{\rm sc,\ell}}\!\-$chain mapping $0$ on $P\,,$ $M_P$ is a 
$(k_\ell)_*\skew4\hat{\hat\P}^{^{\frak X,\ell}}(P)\-$module and the structural $(k_\ell)^*\-$group homomorphism 
$\skew4\hat{\hat\P}^{^{\frak X,\ell}} (P)\to {\rm GL}_{k_\ell} (M_P)$ induces the  reduced
$(k_\ell)_*\skew4\hat{\hat\P}^{^{\frak X,\ell}} \!\!\-$module  
 $\frak m_\ell^{_{\frak X}}\,\colon \skew4\hat{\hat\P}^{^{\frak X,\ell}} \!\to  k_\ell\-\mod\,.$

\smallskip
Otherwise,  choose a minimal element $U$ in $\frak X$ {\it fully normalized\/} in 
$\F$ and set 
$$\frak Y = \frak X - \{\theta(U)\mid \theta\in \F(P,U)\}
\eqno £2.22.1;$$
thus, it follows from our induction hypothesis that  there exists a 
 {\it reduced $(k_\ell)_*\skew4\hat{\hat\P}^{^{\frak Y,\ell}} \!\!\-$module\/}  
$\frak m_\ell^{_{\frak Y}}\,\colon \skew4\hat{\hat\P}^{^{\frak Y,\ell}}\!\to  
k_\ell\-\mod$ such that $X_{\hat\frak q}$ is the class of 
 $\frak m_\ell \big(\hat\frak q (0)\big)$
in $\G_{k_\ell} \big(\skew4\hat{\hat\P}^{^{\frak Y,\ell}} (\hat\frak q)\big)$
for~any  $\skew4\hat{\hat\P}^{^{\frak Y,\ell}} \!\-$ chain $\hat\frak q\,;$ 
let us set $M = \frak m_\ell^{_{\frak Y}}\!(P)\,.$

\smallskip
If $N_\F (U) = \F$ [9,~Proposition~2.16], we also have $N_\P (U) = \P$ [9,~17.5] and then 
$\skew4\hat{\hat\P}^{^{\frak X,\ell}} \!$ actually coincides with the category 
$\T^{^\frak X}_{\skew4\hat{\hat\P}^{^{\frak X,\ell}} \! (U)}$  where $\frak X$ is the set of objects and where,   for a pair of subgroups $Q$ and $R$ in~$\frak X\,,$  the $(k_\ell)^*\-$set of morphisms from $R$ to $Q$ is the 
 {\it $\skew4\hat{\hat\P}^{^{\frak X,\ell}} \!\-$transporter\/} 
 $$\T^{^\frak X}_{\skew4\hat{\hat\P}^{^{\frak X,\ell}} \! (U)} (Q,R) = 
 \{\hat x\in \skew4\hat{\hat\P}^{^{\frak X,\ell}} \! (U)\mid 
 \hat x\.\skew2\hat{\hat\tau}^{_{\rm sc}}_U (R)\.\hat x^{-1}\i 
 \skew2\hat{\hat\tau}^{_{\rm sc}}_U (Q)\}
 \eqno £2.22.2.$$
 Indeed, since $N_\P (U) = \P\,,$ any $\skew4\hat{\hat\P}^{^{\frak X,\ell}} \!\-$morphism 
 $\hat x\,\colon R\to Q$ comes from a $\skew4\hat{\hat\P}^{^{\frak X,\ell}} \!\-$morphism $\hat x^U$ from $R\.U$ to $Q\. U$
stabilizing $U$ and fulfilling (cf.~£2.20)
$$\hat x^U\.\,\hat i^{R\.U}_{R} = \hat i^{Q\.U}_{Q}\.\,\hat x
\eqno £2.22.3;$$
moreover, since $\hat i^{R\.U}_{R}$ is an epimorphism, this equality determines  $\hat x^U$ and the element 
$\frak t_U^{^\frak X} (\hat x)$ of $\skew4\hat{\hat\P}^{^{\frak X,\ell}}\! (U)$ induced by $\hat x^U$ clearly belongs to 
$\T^{^\frak X}_{\skew4\hat{\hat\P}^{^{\frak X,\ell}} (U)} (Q,R)\,;$ thus, this correspondence defines a functor
$$\frak t_U^{^\frak X} : \skew4\hat{\hat\P}^{^{\frak X,\ell}}\!\too 
\T^{^\frak X}_{\skew4\hat{\hat\P}^{^{\frak X,\ell}} (U)}
\eqno £2.22.4\phantom{.}$$
compatible with the structural functors to $\skew4\hat{\hat\F}^{^{\frak X,\ell}}\,,$ and then it is easily checked that this functor is an equivalence of categories.

\smallskip
Now, still denoting by $U$ the {\it $\skew4\hat{\hat\P}^{^{\frak X,\ell}}\!\-$chain\/} mapping $0$ on $U$ and  considering the 
$(k_\ell)_*\skew4\hat{\hat\P}^{^{\frak X,\ell}}(U)\-$module $M_U$ above 
and the structural $(k_\ell)^*\-$group  homomorphism 
$$\hat\rho_U : \skew4\hat{\hat\P}^{^{\frak X,\ell}}\! (U)\too {\rm GL}_{k_\ell} (M_U)
\eqno £2.22.5,$$
it suffices to set $\frak m_\ell^{_{\frak X}} (Q) = M_U = \frak m_\ell^{_{\frak X}} (R)$  and 
$\frak m_\ell^{_{\frak X}} (\hat x) = \hat\rho_U\big( \frak t_U^{^\frak X}(\hat x)\big)\,,$ for any 
$\hat x\in \skew4\hat{\hat\P}^{^{\frak X,\ell}}\! (Q,R)\,,$ to get a  reduced
$(k_\ell)_*\skew4\hat{\hat\P}^{^{\frak X,\ell}}\!\!\-$module  
$$\frak m_\ell^{_{\frak X}} : \skew4\hat{\hat\P}^{^{\frak X,\ell}}\!
\too  k_\ell\-\mod
\eqno £2.22.6;$$   
we claim that it fulfills the announced condition.
\eject

\smallskip
Indeed,  any $\skew4\hat{\hat\P}^{^{\frak X,\ell}}\!\-$chain 
$\hat\frak q\,\colon \Delta_n\to \skew4\hat{\hat\P}^{^{\frak X,\ell}}\!$ can be extended to a 
$\skew4\hat{\hat\P}^{^{\frak X,\ell}}\!\-$chain $\hat\frak q^U\,\colon \Delta_{n+1}\to \skew4\hat{\hat\P}^{^{\frak X,\ell}}\!$ sending $n+1$ to $\hat\frak q (n)\.U$ and $n\bullet n\!+\!1$ to 
$\hat i^{\hat\frak q (n)\.U}_{\hat\frak q (n)} \,,$ and we have an obvious 
$\ch^* (\skew4\hat{\hat\P}^{^{\frak X,\ell}}\!)\-$morphism 
$$({\rm id}_{\hat\frak q},\delta^n_{n+1}) : (\hat\frak q^U,\Delta_{n+1})\too (\hat\frak q,
\Delta_n)
\eqno £2.22.7;$$
similarly, the {\it $\skew4\hat{\hat\P}^{^{\frak X,\ell}}\!\-$chain\/} $U$ can be extended to a $\skew4\hat{\hat\P}^{^{\frak X,\ell}}\!\-$chain
$U^{\hat\frak q (n)}\,\colon \Delta_1\to \skew4\hat{\hat\P}^{^{\frak X,\ell}}\!$
 and we have  obvious $\ch^* (\skew4\hat{\hat\P}^{^{\frak X,\ell}})\-$morphisms
 $$(\hat\frak q^U,\Delta_{n+1})\too (\hat\frak q (n)\.U,\Delta_0)
 \longleftarrow (U^{\hat\frak q (n)},\Delta_1)\too 
 (U,\Delta_0)
 \eqno £2.22.8;$$
 consequently, since we have 
 $$\frak t_U^{^\frak X}\big(\skew4\hat{\hat\P}^{^{\frak X,\ell}}\!
 (\hat\frak q)\big) = \frak t_U^{^\frak X}
 \big(\skew4\hat{\hat\P}^{^{\frak X,\ell}}\!(\hat\frak q^U)\big)\i 
\frak t_U^{^\frak X}\big(\skew4\hat{\hat\P}^{^{\frak X,\ell}}\!
(\hat\frak q (n)\.U)\big) 
 = \frak t_U^{^\frak X}\big(\skew4\hat{\hat\P}^{^{\frak X,\ell}}\!
 (U^{\hat\frak q (n)})\big)
 \eqno £2.22.9\phantom{.}$$
and the family $\{X_{\hat\frak q}\}_{\hat\frak q}$ belongs to
 $\G_{k_\ell} (\skew4\hat{\hat\P}^{^{\frak X,\ell}})\,,$ it follows from our choice of $\ell$ that  we still have the $(k_\ell)_*\skew4\hat{\hat\P}^{^{\frak X,\ell}}\! (\hat\frak q)\-$module isomorphisms (cf.~£2.4.2)
$$\eqalign{M_{\hat\frak q} &\cong M_{\hat\frak q^U}\cong 
 {\rm Res}_{\skew4\hat{\hat\P}^{^{\frak X,\ell}} \!(\hat\frak q)} 
 (M_{\hat\frak q (n)\.U})\cong 
{\rm Res}_{\skew4\hat{\hat\P}^{^{\frak X,\ell}}\!(\hat\frak q)}
(M_{U^{\hat\frak q (n)}})\cr 
&\cong {\rm Res}_{\skew4\hat{\hat\P}^{^{\frak X,\ell}}\!(\hat\frak q)}(M_U) = \frak m_\ell^{_{\frak X}} \!\big(\hat\frak q (0)\big)\cr}
\eqno £2.22.10.$$

\smallskip
Actually,  any reduced  $(k_\ell)_*\skew4\hat{\hat\P}^{^{\frak X,\ell}}\!\-$module 
$\frak n_\ell^{_{\frak X}}\,\colon \skew4\hat{\hat\P}^{^{\frak X,\ell}}\!
\to  k_\ell\-\mod$  fulfilling this condition determines a 
$(k_\ell)_*\skew4\hat{\hat\P}^{^{\frak X,\ell}}\! (U)\-$module structure on 
$\frak n_\ell^{_{\frak X}} (U)$ which, according to our choice of $\ell\,,$  is isomorphic to $M_U\,;$ then, identifying to each other the 
$(k_\ell)_*\skew4\hat{\hat\P}^{^{\frak X,\ell}}\! (U)\-$modules  
$\frak n_\ell^{_{\frak X}} (U)$ and $M_U\,,$ it is easily checked from~£2.22.4 above that  we have the equality $\frak n_\ell^{_{\frak X}} = \frak m_\ell^{_{\frak X}}\,.$
  
\smallskip 
From now on, we assume that $N_\F (U) \not= \F\,;$ then, arguing by induction on the size of $\F\,,$ we may assume that  there exists a unique isomorphism class of reduced $(k_\ell)_*N_{\skew4\hat{\hat\P}^{^{\frak X,\ell}}}\! (U)\-$modules 
$\frak m_{U,\ell}^{_{\frak X}}\,\colon N_{\skew4\hat{\hat\P}^{^{\frak X,\ell}}}\! (U)\to k_\ell\-\mod$ fulfilling the announced  condition; moreover,  $N_{\skew4\hat{\hat\P}^{^{\frak Y,\ell}}}\! (U)$ is a $N_{\F^{^{\frak Y}}} (U)\-$sublocality 
of~$\skew4\hat{\hat\P}^{^{\frak Y,\ell}}\!$ and the restrictions to $N_{\skew4\hat{\hat\P}^{^{\frak Y,\ell}}}\! (U)$ 
of~$\frak m_{U,\ell}^{_{\frak X}}$ and of $\frak m_\ell^{_\frak Y}$ define two reduced
$(k_\ell)_*N_{\skew4\hat{\hat\P}^{^{\frak Y,\ell}}}\! (U)\-$modules,
both fulfilling the announced condition; consequently,  we may assume that 
$\frak m_{U,\ell}^{_{\frak X}}\! \big(N_P(U)\big) = M$ and then, up to a conjugation by a suitable element of~${\rm GL}_{k_\ell} (M)\,,$ that 
$\frak m_{U,\ell}^{_{\frak X}} (\hat y) = \frak m^{_\frak Y} (\hat y)$
for any $N_{\skew4\hat{\hat\P}^{^{\frak Y,\ell}}}\! (U)\-$morphism $\hat y\,\colon R\to Q\,.$

\smallskip
 For any $V\in \frak X -\frak Y$ {\it fully normalized\/} in $\F\,,$  by [9,~Corollary~2.13] there is a $\skew4\hat{\hat\P}^{^{\frak Y,\ell}}\!\-$isomorphism $\hat y\,\colon N_P (U)\to N_P(V)$ fulfilling $\varphi_{\hat y} (U) = V\,,$ where 
 $\varphi_{\hat y}$ is the image of $\hat y$ in $\F\big(N_P (V),N_P (U)\big)\,;$ then, we get a reduced $(k_\ell)_*N_{\skew4\hat{\hat\P}^{^{\frak X,\ell}}}\! (V)\-$mo-dule  $\frak m_{V,\ell}^{_{\frak X}}\,
 \colon N_{\skew4\hat{\hat\P}^{^{\frak X,\ell}}}\! (V)\to k_\ell\-\mod$ sending any $N_{\skew4\hat{\hat\P}^{^{\frak X,\ell}}}\! (V)\-$morphism 
$\hat x\,\colon R\to Q$~to 
$$\frak m_\ell^{_\frak Y}\!(\hat y)^{-1}\circ \frak m_{U,\ell}^{_{\frak X}}\big(\hat r_{Q,\varphi_{\hat y}^{-1}(Q)}^{N_P (V),N_P (U)}(\hat y)^{-1}\. \hat x\. 
\hat r_{R,\varphi_{\hat y}^{-1}(R)}^{N_P (V),N_P (U)}(\hat y)\big)\circ 
\frak m_\ell^{_\frak Y}\!(\hat y)
\eqno £2.22.11\phantom{.}$$
\eject
\noindent
and it is clear that $\frak m_{V,\ell}^{_{\frak X}}\!$ still fulfills the announced condition; moreover, it does
not depend on our choice of $\hat y$ since, for another choice  $\hat y' = 
\hat y\.\hat s\,,$ the element $\hat s$ belongs to 
$\skew4\hat{\hat\P}^{^{\frak X,\ell}}\! \big(N_P(U)\big)_U$ and therefore 
we have $\frak m_{U,\ell}^{_{\frak X}} (\hat s) = 
\frak m_\ell^{_\frak Y} (\hat s)\,;$ similarly, it is easily checked that,  for any 
$N_{\skew4\hat{\hat\P}^{^{\frak Y,\ell}}}\! (V)\-$morphism 
$\hat x\,\colon R\to Q\,,$ we have~$\frak m_{V,\ell}^{_{\frak X}} (\hat x) 
= \frak m_\ell^{_\frak Y} (\hat x)\,.$

\smallskip
At this point, for any  $V,V'\in \frak X -\frak Y$ {\it fully normalized\/} in $\F\,,$  setting $N = N_P (V)$ and $N' = N_P (V')\,,$ 
it follows from~[9,~condition~2.8.2] that any $\skew4\hat{\hat\P}^{^{\frak X,\ell}}\!\-$morphism $\hat x\,\colon V\to V'$ factorizes as $\hat x = \hat r_{V',V}^{N',N} (\hat y)\. \hat s$ for suitable 
$\hat y$ in $\skew4\hat{\hat\P}^{^{\frak X,\ell}}\! (N',N)_{V',V}$ and 
$\hat s$ in $\skew4\hat{\hat\P}^{^{\frak X,\ell}}\! (V)\,;$ then,  in~${\rm GL}_{k_\ell} (M)$ we define 
$$\frak m_\ell^{_{\frak X}} (\hat x) = \frak m_{V,\ell}^{_{\frak X}} (\hat s) \circ 
\frak m_\ell^{_\frak Y} (\hat y)
\eqno £2.22.12;$$
this definition does not depend on our choice since for such another decomposition
$\hat x = \hat r_{V',V}^{N',N} (\skew2\hat{\bar y})\. \skew2\hat{\bar s}\,,$ we get $\skew2\hat{\bar y} = 
\hat y\. \hat t$ and $\skew2\hat{\bar s} = \hat r_{V}^{N} (\hat t)^{-1}\. \hat s$ for a suitable $\hat t$
 in~$\skew4\hat{\hat\P}^{^{\frak X,\ell}}\! (N)_V\,,$ so that we have
 $$\eqalign{\frak m_{V,\ell}^{_{\frak X}} (\skew2\hat{\bar s})\circ 
 \frak m_\ell^{_\frak Y} (\skew2\hat{\bar y})  
&=  \frak m_{V,\ell}^{_{\frak X}} (\hat s) \circ 
\frak m_{V,\ell}^{_{\frak X}} (\hat t)^{-1}
\circ \frak m_\ell^{_{\frak X}}(\hat t)\circ \frak m_\ell^{_\frak Y} (\hat y)\cr
& = \frak m_{V,\ell}^{_{\frak X}} (\hat s) \circ \frak m_\ell^{_\frak Y} (\hat y)\cr}
\eqno £2.22.13.$$

\smallskip
In particular, for any $\hat y\in \skew4\hat{\hat\P}^{^{\frak Y,\ell}}\! (N',N)_{V',V}$ we have
$$\frak m_\ell^{_{\frak X}}\big(\hat r_{V',V}^{N',N} (\hat y)\big) = 
\frak m_\ell^{_{\frak Y}}(\hat y)
\eqno £2.22.14.$$
 More generally, if $Q$ and $Q'$ are a pair of  subgroups of $P$ respectively contained in $N$ and $N'\,,$ and strictly containing $V$ and $V'\,,$ for any 
 $\hat x\in \skew4\hat{\hat\P}^{^{\frak Y,\ell}}\! (Q',Q)_{V',V}$ we 
claim that 
$$\frak m_\ell^{_{\frak X}}\big(\hat r_{V',V}^{Q',Q} (\hat x)\big) = 
\frak m_\ell^{_{\frak Y}}(\hat x)
\eqno £2.22.15;$$
indeed,  it follows from [9,~condition~2.8.2] that $\hat r^{Q',Q}_{V',V} (\hat x) =
\hat r^{N',N}_{V',V} (\hat y)\. \hat z$ for suitable elements
$\hat y\in\skew4\hat{\hat\P}^{^{\frak Y,\ell}}\! (N',N)_{V',V}$ and 
$\hat z\in \skew4\hat{\hat\P}^{^{\frak Y,\ell}}\! (V)\,;$ consequently, setting 
$Q'' = \varphi_{\hat y^{-1}}(Q')\i N\,,$ 
we get
$$\hat z = \hat r^{Q'',Q}_{V,V} \big(\hat r^{N,N'}_{Q'',Q'} 
(\hat y^{-1})\. \hat x\big)
\eqno £2.22.16\phantom{.}$$
and therefore, setting $\hat s = \hat r^{N,N'}_{Q'',Q'} (\hat y^{-1})\. \hat x$ 
which belongs to $\skew4\hat{\hat\P}^{^{\frak Y,\ell}}\! (Q'',Q)_{V,V}\,,$ 
we still get $\hat x = \hat r^{N',N}_{Q',Q''} (\hat y)\. \hat s\,;$ hence, we obtain
$$\hat r_{V',V}^{Q',Q} (\hat x) = \hat r^{N',N}_{V',V} (\hat y)\. 
\hat r_{V,V}^{Q'',Q} (\hat s)\qq \frak m_\ell^{_{\frak Y}}(\hat x) = 
\frak m_\ell^{_{\frak Y}} (\hat s)\circ \frak m_\ell^{_{\frak Y}} (\hat y)
\eqno £2.22.17\phantom{.}$$
and, since $\hat r_{V,V}^{Q'',Q} (\hat s)$ belongs to $\skew4\hat{\hat\P}^{^{\frak X,\ell}}\! (V)$ and  
$\frak m_{V,\ell}^{_{\frak X}}$ is a  {\it reduced $N_{\skew4\hat{\hat\P}^{^{\frak X,\ell}}} \!(V)\-$mo-dule\/}, we finally have (cf.~£2.22.12)
$$\eqalign{\frak m_\ell^{_{\frak X}}\big(\hat r_{V',V}^{Q',Q} (\hat x)\big) 
&= \frak m_{V,\ell}^{_{\frak X}}\! \big(\hat r_{V,V}^{Q'',Q} (\hat s)\big)\circ 
\frak m_\ell^{_{\frak Y}} (\hat y) = \frak m_{V,\ell}^{_{\frak X}}\! (\hat s)\circ 
\frak m_\ell^{_{\frak Y}} (\hat y)\cr
&= \frak m_\ell^{_{\frak Y}}\! (\hat s)\circ 
\frak m_\ell^{_{\frak Y}} (\hat y) = \frak m_\ell^{_{\frak Y}}(\hat x)\cr}
\eqno £2.22.18.$$
\eject

\smallskip
For another  $V''\in \frak X -\frak Y$ {\it fully normalized\/} in $\F\,,$ 
setting $N'' = N_P (V'')$ and considering a 
$\skew4\hat{\hat\P}^{^{\frak X,\ell}}\!\-$morphism $\hat x'\,\colon V'\to V''\,,$ we claim that 
$$\frak m_\ell^{_{\frak X}} (\hat x'\. \hat x) = 
\frak m_\ell^{_{\frak X}} (\hat x)\circ \frak m_\ell^{_{\frak X}} (\hat x')
\eqno £2.22.19;$$
 indeed, assuming that 
 $$\hat x = \hat r_{V',V}^{N',N} (\hat y)\. \hat s\qq
 \hat x' = \hat r_{V'',V'}^{N'',N'} (\hat y')\. \hat s'
 \eqno £2.22.20\phantom{.}$$
 for suitable   $\hat y\in\skew4\hat{\hat\P}^{^{\frak Y,\ell}}\! 
 (N',N)_{V',V}\,,$ $\hat y'\in \skew4\hat{\hat\P}^{^{\frak Y,\ell}}\! 
 (N'',N')_{V'',V'}\,,$ $\hat s\in \skew4\hat{\hat\P}^{^{\frak X,\ell}}\! (V)$ and 
$\hat s'\in \skew4\hat{\hat\P}^{^{\frak X,\ell}}\! (V')\,,$ we get
$$\eqalign{\frak m_\ell^{_{\frak X}} (\hat x)\circ \frak m_\ell^{_{\frak X}} (\hat x') 
&= \frak m_{V,\ell}^{_{\frak X}} (\hat s) \circ \frak m_\ell^{_\frak Y} (\hat y) \circ 
\frak m_{V',\ell}^{_{\frak X}} (\hat s') \circ \frak m_\ell^{_\frak Y} (\hat y')\cr
= \frak m_{V,\ell}^{_{\frak X}} (\hat s) &\circ \big(\frak m_\ell^{_\frak Y} (\hat y) \circ \frak m_{V',\ell}^{_{\frak X}} (\hat s')\circ \frak m_\ell^{_\frak Y} (\hat y)^{-1}\big) \circ \frak m_\ell^{_\frak Y} (\hat y'\. \hat y)\cr
\hat x'\.\hat x &= \hat r_{V'',V'}^{N'',N'} (\hat y')\. \hat s'\.\hat r_{V',V}^{N',N} (\hat y)\. \hat s\cr
&=   \hat r_{V'',V}^{N'',N} (\hat y'\. \hat y)\. \big(\hat r_{V',V}^{N',N} 
(\hat y)^{-1}\. \hat s'\.\hat r_{V',V}^{N',N} (\hat y)\big)\. \hat s\cr}
\eqno £2.22.21.$$

\smallskip
Moreover, it is clear that $\hat y'\. \hat y$ belongs to 
$\skew4\hat{\hat\P}^{^{\frak Y,\ell}}\! (N'',N)_{V'',V}$ and it follows easily from the very definition of $\frak m_{V',\ell}^{_{\frak X}}$ 
in~£2.22.11 above that the element 
$\hat s'' = \hat r_{V',V}^{N',N} (\hat y)^{-1}\.\hat s'\.\hat r_{V',V}^{N',N} 
(\hat y)$ in $\skew4\hat{\hat\P}^{^{\frak X,\ell}}\! (V)$ fulfills
$$ \frak m_{V,\ell}^{_{\frak X}} (\hat s'') = \frak m_\ell^{_\frak Y} (\hat y) \circ 
\frak m_{V',\ell}^{_{\frak X}} (\hat s')\circ \frak m_\ell^{_\frak Y} (\hat y)^{-1}
\eqno £2.22.22;$$
consequently, from~£2.22.21 we obtain
$$\eqalign{\frak m_\ell^{_{\frak X}} (\hat x)\circ \frak m_\ell^{_{\frak X}} (\hat x') 
&= \frak m_{V,\ell}^{_{\frak X}} (\hat s) \circ  \frak m_{V,\ell}^{_{\frak X}} 
(\hat s'') \circ \frak m_\ell^{_\frak Y} (\hat y'\. \hat y)\cr
&= \frak m_{V,\ell}^{_{\frak X}} (\hat s''\. \hat s) \circ 
\frak m_\ell^{_\frak Y} (\hat y'\. \hat y) = 
\frak m_\ell^{_{\frak X}} (\hat x'\. \hat x)\cr}
\eqno £2.22.23,$$
which proves our claim.

\smallskip 
We are ready to consider any pair of subgroups $V$ and $V'$ in $\frak X -\frak Y\,.$ 
 We clearly have  $N = N_P (V)\not= V$ and it follows from [9,~Proposition~2.7] that there is an  $\F\-$morphism $\nu\, \,\colon N\to P$ such that  $\nu (V)$ is {\it fully normalized\/} in~$\F\,;$ moreover, we choose 
 $\hat n\in \skew4\hat{\hat\P}^{^{\frak Y,\ell}}\! 
 \big(\nu (N),N\big)$ lifting the $\F\-$isomorphism $\nu_*$ determined by 
 $\nu\,.$ That is to say, we may assume that
\smallskip
\noindent
£2.22.24\quad {\it There is a pair $(N,\hat n)$ formed by a subgroup $N$ of $P$ which strictly contains and normalizes $V\,,$ and by an element $\hat n$ in 
$\skew4\hat{\hat\P}^{^{\frak Y,\ell}}\!\big(\nu (N),N\big)$ lifting~$\nu_*$ for a  
$\F\-$morphism $\nu\,\colon N\to P$ such that $\nu (V)$ is fully normalized
 in~$\F\,.$\/}
\smallskip
\noindent
 We denote by $\frak N (V)$ the set of such pairs and often we write $\hat n$
instead of $(N,\hat n)\,,$ setting ${}^{\hat n} \!N = \varphi_{\hat n} (N)$
and ${}^{\hat n} \!V = \varphi_{\hat n} (V)\,.$ Then, for any 
$\skew4\hat{\hat\P}^{^{\frak X,\ell}}\!\-$morphism $\hat x\,\colon V\to V'\,,$ we consider 
pairs $(N,\hat n)$ in $\frak N (V)$ and $(N',\hat n')$ in $\frak N (V')$ and, since
${}^{\hat n} \!V$ and ${}^{\hat n'} \!V'$ are both {\it fully normalized\/} in $\F\,,$
we can define
$$\frak m_\ell^{_{\frak X}} \! (\hat x) = \frak m_\ell^{_{\frak Y}} \! (\hat n)\circ 
\frak m_\ell^{_{\frak X}} \!\big(\hat r_{{}^{\hat n'} \!V',V'}^{{}^{\hat n'} \!N',N'} 
(\hat n')\. \hat x\.\hat r_{{}^{\hat n} \!V,V}^{{}^{\hat n} \!N,N} (\hat n)^{-1}\big)\circ \frak m_\ell^{_{\frak Y}} \! (\hat n'^{-1})
\eqno £2.22.25.$$
\eject

\smallskip
This definition is independent of our choices; indeed, for another pair 
$(\bar N,\hat{\bar n})$ in $\frak N (V)\,,$ setting 
$\skew4\bar{\bar N} = \langle N,\bar N\rangle$ and considering a new 
$\F\-$morphism $\psi\, \colon \skew4\bar{\bar N}\to P$  such that $\psi (V)$ is {\it fully normalized\/} in~$\F\,,$  we can obtain a third pair 
$(\skew4\bar{\bar N},\hat m)$ in $\frak N (V)\,;$ then,  
$\hat r_{{}^{\hat m}\! N, N}^{{}^{\hat m} \skew4\bar{\bar N},\skew4\bar{\bar N}}(\hat m)\.\hat n^{-1}$ and $\hat r_{{}^{\hat m}\!\bar N,\bar N}^{{}^{\hat m}\skew4\bar{\bar N},\skew4\bar{\bar N}}(\hat m)\.\hat{\bar n}^{-1}$
respectively
belong to $\skew4\hat{\hat\P}^{^{\frak Y,\ell}}\! 
({}^{\hat m}\! N,{}^{\hat n}\! N)$ and to $\skew4\hat{\hat\P}^{^{\frak Y,\ell}}\!
({}^{\hat m}\!\bar N,{}^{\hat{\bar n}}\! \bar N\big)\,;$
in particular, since ${}^{\hat n}  V\,,$ ${}^{\hat{\bar n}}  V$ and ${}^{\hat m} V$ are {\it fully normalized\/} in~$\F\,,$ we get
$$\eqalign{\frak m_\ell^{_{\frak Y}} \! (\hat m) 
&= \frak m_\ell^{_{\frak Y}} \! \big((\hat r_{{}^{\hat m}\! N, N}^{{}^{\hat m} \skew4\bar{\bar N},\skew4\bar{\bar N}}(\hat m)\.\hat n^{-1})\.\hat n\big)\cr
&=  \frak m_\ell^{_{\frak Y}} \!(\hat n)\circ \frak m_\ell^{_{\frak Y}} \! 
(\hat r_{{}^{\hat m}\! N, N}^{{}^{\hat m} \skew4\bar{\bar N},\skew4\bar{\bar N}}(\hat m)\. \hat n^{-1})\cr
&= \frak m_\ell^{_{\frak Y}} \!(\hat n)\circ \frak m_\ell^{_{\frak X}} \! 
\big(\hat r_{{}^{\hat m}\! V, V}^{{}^{\hat m} 
\skew4\bar{\bar N},\skew4\bar{\bar N}}(\hat m)\.
\hat r_{{}^{\hat n} \!V,V}^{{}^{\hat n} \!N,N} (\hat n^{-1})\big)\cr}
\eqno £2.22.26\phantom{.}$$
and, {\it mutatis mutandis\/}, we still get 
$$\frak m_\ell^{_{\frak Y}} \! (\hat m) = \frak m_\ell^{_{\frak Y}} \!(\hat{\bar n})\circ \frak m_\ell^{_{\frak X}} \! \big(\hat r_{{}^m\! V, V}^{{}^m 
\skew4\bar{\bar N},\skew4\bar{\bar N}}(\hat m)\.\hat r_{{}^{\hat{\bar n}} 
\!V,V}^{{}^{\hat{\bar n}} \! \bar N,\bar N} (\hat{\bar n}^{-1})\big)
\eqno £2.22.27.$$
Consequently, we obtain
$$\frak m_\ell^{_{\frak X}} \! (\hat x) = \frak m_\ell^{_{\frak Y}} \! (\hat{\bar  n})\circ \frak m_\ell^{_{\frak X}} \!
\big(\hat r_{{}^{\hat n'} \!V',V'}^{{}^{\hat n'} \!N',N'} 
(\hat n')\. \hat x\.\hat r_{{}^{\hat{\bar n}} \!V,V}^{{}^{\hat{\bar n}} 
\!\bar N,\bar N} (\hat{\bar n}^{-1})\big)\circ \frak m_\ell^{_{\frak Y}} 
\! (\hat n'^{-1})\eqno £2.22.28.$$
Symmetrically, we can replace  $(\bar N',\hat n')$ for another pair 
$(\bar N',\hat{\bar n}')$ in $\frak N (V')\,.$

\smallskip
Moreover,  equality~£2.22.15 still holds with  this general definition; indeed,
for any pair of subgroups  $Q$ and $Q'$ of $P$ respectively normalizing and strictly containing $V$ and $V'\,,$ we claim that
$$\frak m_\ell^{_{\frak X}}\big(\hat r_{V',V}^{Q',Q} (\hat x)\big) = 
\frak m_\ell^{_{\frak Y}}(\hat x)
\eqno £2.22.29;$$
 indeed, it is clear that we have pairs $(Q,\hat n)$ 
in $\frak N (V)$ and $(Q',\hat n')$ in $\frak N (V')\,,$ and by the very 
definition~£2.22.25 and by  equality~£2.22.15 we have
$$\eqalign{&\frak m_\ell^{_{\frak X}}\big(\hat r_{V',V}^{Q',Q} (\hat x)\big)\cr
& = \frak m_\ell^{_{\frak Y}} \! (\hat n)\circ \frak m_\ell^{_{\frak X}} 
\!\big(\hat r_{{}^{\hat n'} \!V',V'}^{{}^{\hat n'} \!Q',Q'} (\hat n')\. 
\hat r_{V',V}^{Q',Q} (\hat x)\.\hat r_{{}^{\hat n} \!V,V}^{{}^{\hat n} 
\!Q,Q} (\hat n^{-1})\big)\circ \frak m_\ell^{_{\frak Y}} \! (\hat n'^{-1})\cr
&= \frak m_\ell^{_{\frak Y}} \! (\hat n)\circ \frak m_\ell^{_{\frak X}} 
\!\big(\hat r_{{}^{\hat n'} \!V',{}^{\hat n} \!V}^{{}^{\hat n'} \!Q',{}^{\hat n} \!Q} (\hat n'\.  \hat x\.\hat n^{-1})\big)\circ \frak m_\ell^{_{\frak Y}} 
\! (\hat n'^{-1})\cr
&= \frak m_\ell^{_{\frak Y}} \! (\hat n)\circ \frak m_\ell^{_{\frak Y}} 
\! (\hat n'\.  \hat x\.\hat n^{-1})\circ \frak m_\ell^{_{\frak Y}} 
\! (\hat n'^{-1}) = \frak m_\ell^{_{\frak Y}}(\hat x)\cr}
\eqno £2.22.30.$$

\smallskip
Once again, for another  $V''\in \frak X -\frak Y\,,$ setting $N'' = N_P (V'')$ and considering a $\skew4\hat{\hat\P}^{^{\frak X,\ell}}\!\-$morphism 
$\hat x'\,\colon V'\to V''\,,$ we claim that 
$$\frak m_\ell^{_{\frak X}} (\hat x'\. \hat x) = 
\frak m_\ell^{_{\frak X}} (\hat x)\circ \frak m_\ell^{_{\frak X}} (\hat x')
\eqno £2.22.31;$$
indeed, considering a pair $(N'',\hat n'')$ in $\frak N (V'')$ and setting 
$\hat x'' = \hat x'\. \hat x\,,$ from the very 
definition~£2.22.25 we get
$$\eqalign{&\hskip-5pt\frak m_\ell^{_{\frak X}} \! (\hat x) 
= \frak m_\ell^{_{\frak Y}} \! (\hat n)\circ \frak m_\ell^{_{\frak X}} 
\!\big(\hat r_{{}^{\hat n'} \!V',V'}^{{}^{\hat n'} \!N',N'} 
(\hat n')\. \hat x\.\hat r_{{}^{\hat n} \!V,V}^{{}^{\hat n} \!N,N} 
(\hat n^{-1})\big)\circ \frak m_\ell^{_{\frak Y}} \! (\hat n'^{-1})\cr
&\hskip-5pt\frak m_\ell^{_{\frak X}} \! (\hat x') 
= \frak m_\ell^{_{\frak Y}} \! (\hat n')\circ \frak m_\ell^{_{\frak X}} 
\!\big(\hat r_{{}^{\hat n''} \!V'',V''}^{{}^{\hat n''} \!N'',N''} 
(\hat n'')\. \hat x'\.\hat r_{{}^{\hat n'} \!V',V'}^{{}^{\hat n'} \!N',N'} 
(\hat n'^{-1})\big)\circ \frak m_\ell^{_{\frak Y}} \! (\hat n''^{-1})\cr
&\hskip-5pt\frak m_\ell^{_{\frak X}} \! (\hat x'') 
= \frak m_\ell^{_{\frak Y}} \! (\hat n)\circ \frak m_\ell^{_{\frak X}} 
\!\big(\hat r_{{}^{\hat n''} \!V'',V''}^{{}^{\hat n''} \!N'',N''} 
(\hat n'')\. \hat x''\.\hat r_{{}^{\hat n} \!V,V}^{{}^{\hat n} \!N,N} 
(\hat n^{-1})\big)\circ \frak m_\ell^{_{\frak Y}} \! (\hat n''^{-1})\cr}
\eqno £2.22.32\phantom{.}$$
and it follows from equality~£2.22.19 that the composition of the first 
and the second equalities above coincides with the third one.

\smallskip
At this point, we are able to complete the definition of the  reduced $(k_\ell)_*\skew4\hat{\hat\P}^{^{\frak X,\ell}}\!\-$module  
$\frak m_\ell^{_{\frak X}}\,\colon \skew4\hat{\hat\P}^{^{\frak X,\ell}}\!\!\to  k_\ell\-\mod$  fulfilling the announced condition. For any $\skew4\hat{\hat\P}^{^{\frak X,\ell}}\!\!\-$morphism $\hat x\,\colon R\to Q$ either $R$ belongs to $\frak Y$ and we simply set $\frak m_\ell^{_{\frak X}} (\hat x) = \frak m_\ell^{_{\frak Y}} (\hat x)\,,$ or $R$ belongs to $\frak X -\frak Y$ and, denoting by $R_*$ the image of $R$ in $Q$ and by $\hat x_*\,\colon R\cong R_*$ the 
$\skew4\hat{\hat\P}^{^{\frak X,\ell}}\!\-$isomorphism determined by $\hat x\,,$ we set $
\frak m_\ell^{_{\frak X}} (\hat x) = \frak m_\ell^{_{\frak X}} (\hat x_*)$ (cf.~£2.22.25); note that if~$Q$ contains $R$ then we have $\frak m_\ell^{_{\frak X}} \!\big(\hat i^Q_R \big) = {\rm id}_M\,.$ Moreover, we claim that for another $\skew4\hat{\hat\P}^{^{\frak X,\ell}}\!\-$morphism 
$\hat y\,\colon T\to R$ we have
$$\frak m_\ell^{_{\frak X}} \! (\hat x\. \hat y) = 
\frak m_\ell^{_{\frak X}} \! (\hat y)\circ \frak m_\ell^{_{\frak X}} \! (\hat x)
\eqno £2.22.33;$$
indeed, if $T$ belongs to $\frak Y$ then we just have
$$\frak m_\ell^{_{\frak X}} \! (\hat x\. \hat y) = 
\frak m_\ell^{_{\frak Y}} \! (\hat x\. \hat y) = \frak m_\ell^{_{\frak Y}} \! (\hat y)
\circ \frak m_\ell^{_{\frak Y}} \! (\hat x) = \frak m_\ell^{_{\frak X}} \! (\hat y)
\circ \frak m_\ell^{_{\frak X}} \! (\hat x)
\eqno £2.22.34;$$
If $R$ belongs to $\frak X -\frak Y$ then $\hat y$ is a 
$\skew4\hat{\hat\P}^{^{\frak X,\ell}}\!\-$isomorphism and, with the notation above, we have
$T_* = R_*$ and $(\hat x\. \hat y)_* = \hat x_*\.\, \hat y\,;$ in this case,  from equality~£2.22.19 we get
$$\eqalign{\frak m_\ell^{_{\frak X}} \! (\hat x\. \hat y) 
&= \frak m_\ell^{_{\frak X}} \! \big((\hat x\. \hat y)_*\big) = 
\frak m_\ell^{_{\frak X}} \! (\hat x_*\.\, \hat y) =
\frak m_\ell^{_{\frak X}} \! (\hat y) \circ \frak m_\ell^{_{\frak X}} \! (\hat x_*)\cr
& = \frak m_\ell^{_{\frak X}} \! (\hat y)\circ \frak m_\ell^{_{\frak X}} \! (\hat x)\cr}
\eqno £2.22.35.$$

\smallskip
Finally, assume that $T\in \frak X -\frak Y$ and $R\in \frak Y\,,$ denote by $T_*$
and $T_{**}\i R_*$ the respective images of $T$ in $R$ and $Q\,,$ and by $\hat x_{**}\,\colon T_*\to T_{**}$ 
the $\skew4\hat{\hat\P}^{^{\frak X,\ell}}\!\-$iso-morphism fulfilling 
$$\hat x_*\,\.\,\hat i^R_{T_*}  = \hat i^{R_*}_{T_{**}} \,\.\,\hat x_{**}
\eqno £2.22.36;$$ 
then, setting  $\bar R = N_R (T_*)$ and $\bar R_* = N_{R_*} (T_{**})\,,$
it follows from 2.23.16 and £2.22.31 that we have
$$\eqalign{\frak m_\ell^{_{\frak X}} \! (\hat x\. \hat y) 
&= \frak m_\ell^{_{\frak X}} \! (\hat x_{**}\,\.\, \hat y_*) =  \frak m_\ell^{_{\frak X}} \! (\hat y_*)\circ \frak m_\ell^{_{\frak X}} \! (\hat x_{**})\cr
& = \frak m_\ell^{_{\frak X}} \! (\hat y)\circ \frak m_\ell^{_{\frak X}} \! 
\Big(\hat r_{T_{**},T_*}^{\bar R_*,\bar R} \big(\hat r_{\bar R_*,\bar R}^{Q,R}
(\hat x)\big)\Big)\cr
&= \frak m_\ell^{_{\frak X}} \! (\hat y)\circ \frak m_\ell^{_{\frak Y}} \! 
 \big(\hat r_{\bar R_*,\bar R}^{Q,R}(\hat x)\big) = 
 \frak m_\ell^{_{\frak X}} \! (\hat y)\circ \frak m_\ell^{_{\frak Y}} \!  (\hat x)\cr}
\eqno £2.22.37.$$
It is not difficult to check that the functor $\frak m_\ell^{_{\frak X}}\,\colon \skew4\hat{\hat\P}^{^{\frak X,\ell}}\!\to  k_\ell\-\mod$ is a reduced $(k_\ell)_*\skew4\hat{\hat\P}^{^{\frak X,\ell}}\!\-$module  which
 fulfills the announced condition. We are done.
 \eject

 \bigskip
 \noindent
 {\bf Corollary~£2.23.} {\it  The homomorphism $\cat_\K\,\colon
\G (\K_*\hat\P^{^{\rm sc}}\!\-\mod)\to \G_\K (\hat\P^{^{\rm sc}})$ is bijective.\/}

\medskip
\noindent
{\bf Proof:} According to Corollary~£2.17, it suffices to prove the surjectivity. Let $\{X_{\hat\frak q}\}_{\hat\frak q}$ 
where $\hat\frak q$ runs over the set of {\it $\hat\P^{^{\rm sc}}\-$chains\/} be a family belonging 
to~$\G_\K (\hat\P^{^{\rm sc}})\,;$ according to Remark~£1.11, we actually may assume that $X_{\hat\frak q}$ is the class of a 
$\K_*\hat\P^{^{\rm sc}} (\hat\frak q)\-$module $M_{\hat\frak q}$ for any {\it $\hat\P^{^{\rm sc}}\!\-$chain\/} 
$\hat\frak q\,\colon \Delta_n\to \hat\P^{^{\rm sc}}$ and then, since the isomorphisms~£2.18.1,~£2.19.1 and~£2.19.2
all preserve the classes of modules, it follows from isomorphism~£2.20.3 that $\{X_{\hat\frak q}\}_{\hat\frak q}$ determines a family  $\{X_{\hat\frak q^{^\ell}}\}_{\hat\frak q^{^\ell}}$ in~$\G_{k_\ell} \big(\skew4\hat{\hat\P}^{^{\rm sc,\ell}})$ where  
$\hat\frak q^{^\ell}$ runs over the set of  {\it $\skew4\hat{\hat\P}^{^{\rm sc,\ell}}\!\-$chains\/} and, for such a 
$\hat\frak q^{^\ell}\,,$ $X_{\hat\frak q^{^\ell}}$ is the class of a 
$(k_\ell)_*\skew4\hat{\hat\P}^{^{\rm sc,\ell}}\! (\hat\frak q)\-$module 
$M_{\hat\frak q^{^\ell}}\,.$ Consequently, it follows from Theorem~£2.22 above that
there exists a reduced $(k_\ell)_*\skew4\hat{\hat\P}^{^{\rm sc,\ell}}\!\-$module
$\frak m_\ell\,\colon \skew4\hat{\hat\P}^{^{\rm sc,\ell}}\!\to k_\ell\-\mod$  
such that $X_{\hat\frak q^{^\ell}}$ is the class of 
$\frak m_\ell \big(\hat\frak q^{^\ell} (0)\big)$ in 
$\G_{k_\ell} \big(\skew4\hat{\hat\P}^{^{\rm sc,\ell}} (\hat\frak q^{^\ell})\big)$
for~any  $\skew4\hat{\hat\P}^{^{\rm sc,\ell}} \!\-$chain $\hat\frak q^{^\ell}\,.$

\smallskip
In particular, setting $M_\ell = \frak m_\ell (P)$ and denoting by  $\hat G (\frak m_\ell)$ the $(k_\ell)^*\-$sub-group of 
the {\it finite\/} group ${\rm GL}_{k_\ell} (M_\ell)$ generated  by~$\frak m_\ell (\hat x)$ where  $\hat x\,\colon R\to Q$ runs over the set of  $\skew4\hat{\hat\P}^{^{\rm sc,\ell}}\!\-$morphisms, $M_\ell$ becomes 
a $(k_\ell)_*\hat G (\frak m_\ell)\-$module, determining an element
$X_\ell$ in $\G_{k_\ell} \big(\hat G (\frak m_\ell)\big)\,;$ since the {\it Brauer decomposition map\/}
$$\G_{\K_\ell} \big(\hat G (\frak m_\ell)\big)\too \G_{k_\ell} \big(\hat G (\frak m_\ell)\big)
\eqno £2.23.1\phantom{.}$$
is surjective, $X_\ell$ can be lifted to some element $\hat X_\ell\in \G_{\K_\ell} \big(\hat G (\frak m_\ell)\big)$ which is then the difference of the classes of suitable $(\K_\ell)_*\hat G (\frak m_\ell)\-$modules $\hat M'_\ell$ and~$\hat M''_\ell\,.$ 

\smallskip
But, since we have an obvious $(k_\ell)^*\-$functor from the $(k_\ell)^*\-$category $\skew4\hat{\hat\P}^{^{\rm sc,\ell}} \!$
to the $(k_\ell)^*\-$category with a unique object $\emptyset$ and $(k_\ell)^*\-$group of automorphisms 
$\hat G (\frak m_\ell)\,,$ $\hat M'_\ell$ and~$\hat M''_\ell$ determine respective 
{\it  reduced $(\K_\ell)_*\skew4\hat{\hat\P}^{^{\rm sc,\ell}}\!\-$modules\/} 
$\hat\frak m'_\ell\,\colon \skew4\hat{\hat\P}^{^{\rm sc,\ell}}\!\to \K_\ell\-\mod$ 
and $\hat\frak m''_\ell\,\colon \skew4\hat{\hat\P}^{^{\rm sc,\ell}}\!\to \K_\ell\-\mod\,;$ then, since 
$\skew4\hat{\hat\P}^{^{\rm sc}}\i \skew4\hat{\hat\P}^{^{\rm sc,\ell}}$ and $\K_\ell\i \skew4\hat{\hat\K}$ (cf.~£2.19),
$\hat\frak m'_\ell$ and $\hat\frak m''_\ell$  determine respective 
{\it  reduced $\skew4\hat{\hat\K}_*\skew4\hat{\hat\P}^{^{\rm sc}}\!\-$modules\/} 
$\hat{\hat\frak m}'\,\colon \skew4\hat{\hat\P}^{^{\rm sc}}\!\to \skew4\hat{\hat\K}\-\mod$ 
and $\hat{\hat\frak m}''\,\colon \skew4\hat{\hat\P}^{^{\rm sc}}\!\to \skew4\hat{\hat\K}\-\mod\,;$
finally, since we actually may assume that $\skew4\hat{\hat\K}$ contains $\K$ and that $\K$ is big enough,
$\hat{\hat\frak m}'$ and $\hat{\hat\frak m}''$ come from respective {\it  reduced $\K_*\hat\P^{^{\rm sc}}\!\-$modules\/} 
$\frak m'\,\colon \hat\P^{^{\rm sc}}\!\to \K\-\mod$  and $\frak m''\,\colon \hat\P^{^{\rm sc}}\!\to \K\-\mod\,,$
and it is easily checked that the diference between their images 
in~$\G_\K (\hat\P^{^{\rm sc}})$ coincides with the starting family $\{X_{\hat\frak q}\}_{\hat\frak q}\,.$

  \bigskip
 \noindent
 {\bf Corollary~£2.24.} {\it There exists a  $\K_*\hat\P^{^{\rm sc}}\-$module 
 $\frak m\,\colon \hat\P^{^{\rm sc}}\! \to \K\-\mod$ such that, for any 
 $\hat\P^{^{\rm sc}}\!\-$chain $\hat\frak q\,\colon \Delta_n\to 
 \hat\P^{^{\rm sc}}\,,$ the class of $\frak m \big(\hat\frak q (0)\big)$
 in $\G_\K \big(\hat\P^{^{\rm sc}}(\hat\frak q)\big)$ is a multiple of the regular
 $\K_*\hat\P^{^{\rm sc}}\! (\hat\frak q)\-$module. In particular, the $k^*\-$functor from $\hat\P^{^{\rm sc}}\!$ to the 
 $k^*\-$category over one object with automorphism $k^*\-$group  $\hat G (\frak m)$ is faithful.\/}
 \eject
 
\medskip
\noindent
{\bf Proof:}  From Remark~£1.11 we already know that,  denoting by~$R_{\hat\frak q}$ the class in 
$\G_{k_\ell} \big(\skew4\hat{\hat\P}^{^{\rm sc,\ell}}\! (\hat\frak q)\big)$ of the {\it regular\/} 
 $(k_\ell)_*\skew4\hat{\hat\P}^{^{\rm sc,\ell}}\! (\hat\frak q)\-$module 
 $(k_\ell)_*\skew4\hat{\hat\P}^{^{\rm sc,\ell}}\! (\hat\frak q)\,,$ and choosing
a multiple $m$ of all the orders $\vert\P^{^{\rm sc}}\! (\hat\frak q)\vert$ where $\hat\frak q$ runs over the set of  
{\it $\skew4\hat{\hat\P}^{^{\rm sc,\ell}}\!\-$chains,\/} it is easily checked that  the family $R = \Big\{\displaystyle{m\over \vert\P^{^{\rm sc}}\! (\hat\frak q)\vert}\. R_{\hat\frak q}\Big\}_{\hat\frak q}$ belongs to
$\G_{k_\ell} (\skew4\hat{\hat\P}^{^{\rm sc,\ell}})\,;$
hence, it follows from Theorem~£2.22 that  there exists a reduced $(k_\ell)_*\skew4\hat{\hat\P}^{^{\rm sc,\ell}}\!\-$module
$\frak n_\ell\,\colon \skew4\hat{\hat\P}^{^{\rm sc,\ell}}\!\to k_\ell\-\mod$  
such that $\displaystyle{m\over \vert\P^{^{\rm sc}}\! (\hat\frak q)\vert}\. R_{\hat\frak q}$ is the class of 
$\frak n_\ell \big(\hat\frak q (0)\big)$
in $\G_{k_\ell} \big(\skew4\hat{\hat\P}^{^{\rm sc,\ell}} (\hat\frak q)\big)$
for~any  $\skew4\hat{\hat\P}^{^{\rm sc,\ell}}\! \!\-$ chain $\hat\frak q\,.$

\smallskip
As above, setting $N_\ell = \frak n_\ell (P)$ and denoting by $\hat G (\frak n_\ell)$ the $(k_\ell)^*\-$subgroup 
of the {\it finite\/} group ${\rm GL}_{k_\ell} (N_\ell)$ generated  by~$\frak n_\ell (\hat x)$ where  
$\hat x\,\colon R\to Q$ runs over the set of  $\skew4\hat{\hat\P}^{^{\rm sc,\ell}}\!\-$morphisms, $N_\ell$ becomes 
a $(k_\ell)_*\hat G (\frak n_\ell)\-$module determining an element in $\G_{k_\ell} \big(\hat G (\frak n_\ell)\big)\,;$ 
then, we know that a multiple of this element comes, {\it via\/}  the {\it Brauer decomposition map\/}, from a {\it true\/} $(\K_\ell)_*\hat G (\frak n_\ell)\-$module
and, {\it via\/}  a field $\hat\K$ containing $\K_\ell$ and $\K\,,$  it comes from a $\K_*\hat G (\frak n_\ell)\-$module
$M\,;$ finally, the reduced $\K_*\hat\P^{^{\rm sc}}\-$module  $\frak m\,\colon \hat\P^{^{\rm sc}}\! \to \K\-\mod$ determined by 
$M$ fulfills the announced condition.

\medskip
£2.25.  Let $\frak s_{\hat\P^{^{\rm sc}}} \,\colon \hat\P^{^{\rm sc}}\!\!\to \K\-\mod$ be the 
 direct sum of a set of reduced representatives for the set of isomorphism classes of simple $\K_*\hat\P^{^{\rm sc}}\!\-$modules and denote by~$\hat G(\hat\P^{^{\rm sc}})$ the $k^*\-$subgroup  of  ${\rm GL}_\K \big(\frak s_{\hat\P^{^{\rm sc}}} (P)\big)$ generated  by  $\frak s_{\hat\P^{^{\rm sc}}} (\hat x)$ where  $\hat x\,\colon R\to Q$ runs over the set of   $\hat\P^{^{\rm sc}}\!\-$morphisms and by $\frak f_{\hat\P^{^{\rm sc}}}$ the functor   determined by 
 $\frak s_{\hat\P^{^{\rm sc}}}$ from $\hat\P^{^{\rm sc}}\!\!$ to the   $k^*\-$category over one object with the automorphism 
 $k^*\-$group~$\hat G (\hat\P^{^{\rm sc}})\,,$ which is {\it faithful\/} by Corollary~£2.24 above.

 \bigskip
 \noindent
 {\bf Corollary~£2.26.} {\it With the notation above,  the  functor $\frak f_{\hat\P^{^{\rm sc}}}$ induces an~equivalence of categories from  $\K_*\hat G(\hat\P^{^{\rm sc}})\-\mod$  to~$\K_*\hat\P^{^{\rm sc}}\!\-\mod\,.$ Moreover, the regular representation  of~$\hat G(\hat\P^{^{\rm sc}})$ induces a  $\K_*\hat\P^{^{\rm sc}}\!\-$module 
 $\frak r_{\hat\P^{^{\rm sc}}} \,\colon \hat\P^{^{\rm sc}}\! \to \K\-\mod$ such that, for any  $\hat\P^{^{\rm sc}}\!\-$chain 
 $\hat\frak q\,\colon \Delta_n\to \hat\P^{^{\rm sc}}\,,$ the class of $\frak r_{\hat\P^{^{\rm sc}}} \big(\hat\frak q (0)\big)$
 in $\G_\K \big(\hat\P^{^{\rm sc}}(\hat\frak q)\big)$ is a multiple of the regular
 $\K_*\hat\P^{^{\rm sc}}\! (\hat\frak q)\-$module.\/}
 
 \medskip
 \noindent
 {\bf Proof:} It is clear that,  for any simple $\K_*\hat\P^{^{\rm sc}}\!\-$module 
 $\frak s \,\colon \hat\P^{^{\rm sc}}\!\!\to \K\-\mod\,,$ the $k^*\-$group $\hat G(\hat\P^{^{\rm sc}})$ acts on $\frak s (P)$
 and then $\frak s (P)$ becomes a simple $\K_*\hat G(\hat\P^{^{\rm sc}})\-$ module; more generally, the restriction from 
 $\hat G(\hat\P^{^{\rm sc}})$  to $\hat\P^{^{\rm sc}}\!$ {\it via\/} $\frak f_{\hat\P^{^{\rm sc}}}$ clearly  determines a functor  from 
 $\K_*\hat G(\hat\P^{^{\rm sc}})\-\mod$  to $\K_*\hat\P^{^{\rm sc}}\!\-\mod$ which induces a bijection between the sets of isomorphism classes of simple $\K_*\hat G(\hat\P^{^{\rm sc}})\-$ and $\K_*\hat\P^{^{\rm sc}}\!\-$modules; since both categories are semisimple, this functor is an equivalence of categories. Moreover, since the functor $\frak f_{\hat\P^{^{\rm sc}}}$ 
 is {\it faithful\/},  for any  $\hat\P^{^{\rm sc}}\!\-$chain  $\hat\frak q\,\colon \Delta_n\to \hat\P^{^{\rm sc}}$
 the $k^*\-$group $\hat\P^{^{\rm sc}}\! (\hat\frak q)$ is $k^*\-$isomorphic to a $k^*\-$subgroup 
 of~$\hat G(\hat\P^{^{\rm sc}})\,.$

\bigskip
\noindent
{\bf £3. Categorization for the characteristic $p$ case}
\medskip
£3.1. With the notation in~£2.1 above, this time we set (cf.~£1.10.3)
$$\G_k (\hat\P^{^{\rm sc}}) = 
\lim_{\longleftarrow}\,(\frak g_k\circ \aut_{\hat\P^{^{\rm sc}}})
\eqno £3.1.1;$$
that is to say, $\G_k (\hat\P^{^{\rm sc}}) $ is the subset of elements
$$\{Z_{(\hat\frak q,\Delta_n)}\}_{(\hat\frak q,\Delta_n)}\in \prod_{n\in \Bbb N}\,
\prod_{\hat\frak q\in \Fct (\Delta_n,\hat\P^{^{\rm sc}})} \G_k \big(\hat\P^{^{\rm sc}}\! (\hat\frak q)\big)
\eqno £3.1.2\phantom{.}$$
such that, for any $\hat\P^{^{\rm sc}}\-$morphism $(\nu,\delta)\,\colon 
(\hat\frak q,\Delta_n)\to (\hat\frak r,\Delta_m)\,,$ they fulfill
$${\rm res}_{\aut_{\hat\P^{^{\rm sc}}}(\nu,\delta)}(Z_{(\hat\frak r,\Delta_m)})
= Z_{(\hat\frak q,\Delta_n)}
\eqno £3.1.3.$$
Our purpose in this section is to show that  $\G_k (\hat\P^{^{\rm sc}})$ is the extension to $\O$ of the very {\it Grothendieck group\/} of the category of $k_*\hat G(\hat\P^{^{\rm sc}})\-$modules for the $k^*\-$group~$\hat G (\hat\P^{^{\rm sc}})$
introduced in~£2.25 above.

 \medskip
 £3.2. As in~£2.3 above, we call {\it contravariant $k^*\-$functor\/} 
 $\frak m\,\colon \hat\P^{^{\rm sc}}\!\to k\-\mod$ any functor such that   
  $\frak m[\lambda\.\hat x) = \lambda\.\frak m (\hat x)$
 for any $\hat\P^{^{\rm sc}}\!\-$morphism $\hat x\,\colon R\to Q$ and any 
 $\lambda\in k^*\,;$ once again, any {\it contravariant $k^*\-$functor\/} 
 $\frak m\,\colon \hat\P^{^{\rm sc}}\!\to k\-\mod$ determines a new 
 {\it contravariant $k^*\-$functor\/}  [9,~A3.7.3]
 $$\frak m^{\ch} = \frak m\circ \frak v_{\hat\P^{^{\rm sc}}}: 
 \ch^* (\hat\P^{^{\rm sc}})\too k\-\mod
 \eqno £3.2.1\phantom{.}$$
 sending any {\it $\hat\P^{^{\rm sc}}\!\-$chain\/} 
 $\hat\frak q\,\colon \Delta_n\to \hat\P^{^{\rm sc}}\!$ to 
 $\frak m\big(\hat\frak q (0)\big)$ and any 
 $\ch^*(\hat\P^{^{\rm sc}})\-$morphism
 $$(\hat x,\delta) : (\hat\frak r,\Delta_m)\too (\hat\frak q,\Delta_n)
 \eqno £3.2.2,$$
 where $\delta\,\colon \Delta_n\to \Delta_m$ is an order-preserving map
and $\hat x\,\colon \hat\frak r\circ\delta\cong \hat\frak q$ a 
{\it natural isomorphism\/}, to the $k\-$linear map 
$$\frak m \big(\hat x_0\circ \hat r (0\bullet \delta (0)\big)
:  \frak m\big(\hat\frak q (0)\big)\too \frak m\big(\hat\frak r (0)\big)
\eqno £3.2.3.$$

\medskip
£3.3. As in £2.4 above, in the case where $n = m$ and $\delta$ is the identity map, we get a $k^*\-$compatible 
action of the $k^*\-$group $\aut_{\hat\P^{^{\rm sc}}\!} (\hat\frak q) = \hat\P^{^{\rm sc}}\! (\hat\frak q)$ 
over $\frak m\big(\hat\frak q (0)\big)\,,$ so that $\frak m\big(\hat\frak q (0)\big)$ becomes a 
$k_*\hat\P^{^{\rm sc}}\! (\hat\frak q)\-$module; similarly, 
$\frak m\big(\hat\frak r (0)\big)$ becomes a $k_*\hat\P^{^{\rm sc}}\! (\hat\frak r)\-$module and, 
{\it via\/} the $k^*\-$group homomorphism $\aut_{\hat\P^{^{\rm sc}}\!} (\hat x,\delta)\,,$
$\frak m\big(\hat\frak q (0)\big)$ also becomes the $k_*\hat\P^{^{\rm sc}}\! (\hat\frak r)\-$module 
${\rm Res}_{\aut_{\hat\P^{^{\rm sc}}\!} (\hat x,\delta)}\Big(\frak m\big(\hat\frak q (0)\big)\Big)\,;$ 
then the  $k\-$linear map~£3.2.3  is clearly a  $k_*\hat\P^{^{\rm sc}}\! (\hat\frak r)\-$module homomorphism
 $$\frak m^{\ch} (\hat x,\delta) :  
 {\rm Res}_{\aut_{\hat\P^{^{\rm sc}}\!} (\hat x,\delta)}
\Big(\frak m\big(\hat\frak q (0)\big)\Big)\too \frak m\big(\hat\frak r (0)\big)
\eqno £3.3.1.$$
 Similarly, any {\it natural map\/} $\mu\,\colon \frak m\to \frak m'$ between {\it contravariant $k^*\-$ functors\/}  $\frak m$ and $\frak m'$ from 
 $\hat\P^{^{\rm sc}}\!$ to $k\-\mod$ determines  a new  {\it natural map\/} 
 $$\mu^{\ch} = \mu *\frak v_{\hat\P^{^{\rm sc}}} : 
 \frak m^{\ch}\too \frak m'^{\ch}
 \eqno £3.3.2\phantom{.}$$
 sending any {\it $\hat\P^{^{\rm sc}}\!\-$chain\/} 
 $\hat\frak q\,\colon \Delta_n\to \hat\P^{^{\rm sc}}\!$ to the $k\-$linear map $\mu_{\hat\frak q (0)}$
and it follows from the naturalness of $\mu$ that this map is actually
a $k_*\hat\P^{^{\rm sc}}\! (\hat\frak q)\-$module homomorphism.
\eject

\medskip
£3.4. Once again, we are interested in the {\it contravariant $k^*\-$functors\/} 
 $\frak m\,\colon \hat\P^{^{\rm sc}}\!\to k\-\mod$ --- called {\it reversible\/} --- mapping {\it any\/}  
 $\hat\P^{^{\rm sc}}\!\-$morphism $\hat x\,\colon R\to Q$ on a $k\-$linear isomorphism 
 $\frak m(\hat x)\,\colon \frak m(Q)\cong \frak m(R)\,;$
 in this case, it is quite clear that the {\it contravariant $k^*\-$functor\/}  
 $\frak m^{\ch}$ also maps any $\ch^* (\hat\P^{^{\rm sc}})\-$morphism
 on a $k\-$linear isomorphism and, as in~£2.6.2 above, the converse is also true.
   Moreover, if $\frak m\,\colon \hat\P^{^{\rm sc}}\!\to k\-\mod$ is a
{\it reversible contravariant $k^*\-$ functor\/},  homomorphism~£3.3.1 becomes an isomorphism
 and therefore the restriction map
 $$ \frak g_k \big(\aut_{\hat\P^{^{\rm sc}}\!} (\hat x,\delta)\big) :
 \G_k (\hat\P^{^{\rm sc}}\! (\hat\frak q)\big)\too \G_k (\hat\P^{^{\rm sc}}\! (\hat\frak r)\big)
 \eqno £3.4.1\phantom{.}$$
  sends  the class $Z_{\frak m (\hat\frak q (0))}$ in $\G_k (\hat\P^{^{\rm sc}}\! (\hat\frak q)\big)$ of the 
 $k_*\hat\P^{^{\rm sc}}\! (\hat\frak q)\-$module $\frak m \big(\hat\frak q (0)\big)$ to the class 
 $Z_{\frak m (\hat\frak r (0))}$  in $\G_k (\hat\P^{^{\rm sc}}\! (\hat\frak r)\big)$ of the 
 $k_*\hat\P^{^{\rm sc}}\! (\hat\frak r)\-$module 
 $\frak m \big(\hat\frak r (0)\big)\,.$ That is to say,\break
 the family
 $$\big\{Z_{\frak m (\hat\frak q (0))}\big\}_{(\hat\frak q,\Delta_n)}
 \in \prod_{n\in\Bbb N}\,\prod_{\hat\frak q\in \Fct (\Delta_n,\hat\P^{^{\rm sc}})}
 \G_k (\hat\P^{^{\rm sc}}\! (\hat\frak q)\big) 
 \eqno £3.4.2\phantom{.}$$
 fulfills condition~£3.1.3 and therefore it belongs to $\G_k (\hat\P^{^{\rm sc}})\,;$  let us denote by $Z_\frak m$ this family which, clearly, only depends on the isomorphism class of~$\frak m\,.$

\medskip
£3.5. This time we call {\it $k_*\hat\P^{^{\rm sc}}\!\-$module\/} $\frak m\,\colon \hat\P^{^{\rm sc}}\!\to k\-\mod$ just any {\it contra-variant\/} functor obtained by restriction from a  $k_*\hat G(\hat\P^{^{\rm sc}})\-$module {\it via\/} the functor
$\frak f_{\hat\P^{^{\rm sc}}}$ introduced in~£2.15; thus, denoting by $k_*\hat\P^{^{\rm sc}}\!\-\mod$ the category 
formed by these {\it $k_*\hat\P^{^{\rm sc}}\!\-$modules\/} and by the {\it natural maps\/} between them, we clearly get 
a {\it faithful\/} functor
$$k_*\hat G(\hat\P^{^{\rm sc}})\-\mod\too k_*\hat\P^{^{\rm sc}}\!\-\mod
\eqno £3.5.1\,;$$
moreover, it is easily checked that any {\it natural map\/} $\mu\,\colon \frak m\to \frak m'$ between 
$k_*\hat\P^{^{\rm sc}}\!\-$modules $\frak m$ and $\frak m'$ induces a $k_*\hat G(\hat\P^{^{\rm sc}})\-$morphism
$\mu_P\,\colon \frak m (P)\to \frak m'(P)\,;$ consequently, the functor above is actually an equivalence of categories.
Thus, denoting by $\G (k_*\hat\P^{^{\rm sc}}\!\-\mod)$ the extension to $\O$ of the Grothendieck group 
of~$k_*\hat\P^{^{\rm sc}}\!\-\mod\,,$ the equivalence~£3.5.1 induces an $\O\-$module isomorphism
$$\G (k_*\hat\P^{^{\rm sc}}\!\-\mod)\cong \G_k \big(\hat G(\hat\P^{^{\rm sc}})\big)
\eqno £3.5.2.$$

\medskip
£3.6. Moreover, as above the correspondence sending any {\it $k_*\hat\P^{^{\rm sc}}\!\-$module\/} 
 $\frak m$ to  $Z_\frak m$ induces an $\O\-$module homomorphism
 $$\cat_k : \G (k_*\hat\P^{^{\rm sc}}\!\-\mod)\too \G_k (\hat\P^{^{\rm sc}})
 \eqno £3.6.1\phantom{.}$$
 and it is quite clear that  suitable {\it Brauer decomposition maps\/} [10,~3.4.2] yield the following commutative diagram (cf.~Corollary~£2.23)
 $$\matrix{\G_\K \big(\hat G(\hat\P^{^{\rm sc}})\big)&\cong &\G (\K_*\hat\P^{^{\rm sc}}\!\-\mod) &\buildrel \cat_\K\over\cong &\G_\K (\hat\P^{^{\rm sc}}) \cr
\hskip-40pt{\scriptstyle \delta_{\hat G(\hat\P^{^{\rm sc}})}}\hskip5pt\big\downarrow &&
&&\big\downarrow\hskip5pt{\scriptstyle \delta_{\hat\P^{^{\rm sc}}}}\hskip-25pt\cr
\G_k \big(\hat G(\hat\P^{^{\rm sc}})\big)&\cong &\G (k_*\hat\P^{^{\rm sc}}\!\-\mod) &\buildrel \cat_k\over\too 
 &\G_k (\hat\P^{^{\rm sc}})  \cr}
 \eqno £3.6.2\,.$$
But we already know that, for any finite group $H\,,$ the  {\it Brauer decomposition map\/} 
$\delta_H\,\colon \G_\K (H)\to \G_k (H)$ is surjective and even admits a {\it natural section\/} --- namely, extending any 
{\it Brauer character\/} $\varphi$ of $H$ to the central function over~$H$ mapping $y\in H$ on $\varphi (y_{p'})\,;$ moreover, it is easy to check that all this still holds for $k^*\-$groups. Consequently, the vertical homomorphisms in the diagram above admit sections and, in particular, they are surjective.

\medskip
£3.7. More generally, for any finite group $H$ and any $p\-$element $v$ in~$H\,,$ assuming that $\K$ is big enough for $H$ we have the {\it Brauer general decomposition map\/}
$$\delta^v_H : \G_\K (H)\too \G_k \big(C_H (v)\big)
\eqno £3.7.1\phantom{.}$$
which, following [2, Appendice], can be defined as the composition 
$$\delta^v_H = \delta_{C_H (v)}\circ \omega^v_{C_H (v)}\circ {\rm Res^H_{C_H (v)}}
\eqno £3.7.2\phantom{.}$$
where we consider the {\it $v\-$twist\/}
$$\omega^v_{C_H (v)} : \G_\K \big(C_H (v)\big)\too \G_\K \big(C_H (v)\big)
\eqno £3.7.3\phantom{.}$$ 
determined by the {\it $v\-$translation\/} map  induced by the multiplication by $v$ in the $\O\-$valued functions 
$\F\! ct \big(C_H (v),\O\big)$ [10,~9.2]; then, recall that we have  [2, Appendice]
$$\K\otimes_\O \G_\K (H)\cong \prod_v  \K\otimes_\O \G_k \big(C_H (v)\big)
\eqno £3.7.4\phantom{.}$$
where $v$ runs over a set of representatives for the set of conjugacy classes of $p\-$elements in $H\;.$

\medskip
£3.8. In order to get a  similar result for a $k^*\-$group $\hat H$ with finite $k^*\-$quo-tient $H\,,$ we have to replace {\it $p\-$element\/} by 
{\it local element\/}; we say that a $p\-$element $v$ of $\hat H$ is {\it local\/} if either $v = 1$ or the {\it relative trace map\/}
$${\rm tr}^{\langle v\rangle}_{\langle v^p\rangle} : (k_*\hat H)^{\langle v^p\rangle} \too (k_*\hat H)^{\langle v\rangle}
\eqno £3.8.1\phantom{.}$$
is {\it not\/} surjective. Then, from the standard results on $k^*\-$groups [7,~Proposition~5.15], it is easily checked from isomorphism~3.7.4 that
$$\K\otimes_\O \G_\K (\hat H)\cong \prod_v  \K\otimes_\O \G_k \big(C_{\hat H} (v)\big)
\eqno £3.8.2\phantom{.}$$
where $v$ runs over a set of representatives for the set of conjugacy classes of {\it local elements\/} in $\hat H\,.$

\medskip
£3.9. The point is that in [10,~Theorem~9.3] we prove an analogous result for $\G_\K (\hat\P^{^{\rm sc}})\,;$
explicitly,  choose a set of representatives $\U\i P$ for the set of $\F\-$isomorphism classes of the elements of $P$ in such a way that, for any $u\in \P\,,$ the subgroup $\langle u\rangle$ is fully centralized in $\F$ [9,~Proposition~2.7].\break
\eject
\noindent
For any $u\in \U\,,$ we have the Frobenius $C_P (u)\-$category $C_\F (u)$ [9,~Proposition~2.16] and,  since a 
$C_\F (u)\-$selfcentralizing subgroup of~$C_P (u)$ contains~$u\,,$ so that it is also an $\F\-$selfcentralizing subgroup of~$P\,,$  we write $C_{\F^{^{\rm sc}}} (u)$ instead of $C_{\F} (u)^{^{\rm sc}}\,.$ Then, with obvious notation, it is easily checked that 
$C_{\P} (u)$ and $C_{\P^{^{\rm sc}}} (u)$ are the respective {\it perfect $C_\F (u)\-$ and $C_{\F^{^{\rm sc}}} (u)\-$localities\/}, and we clearly have the {\it finite folded Frobenius $C_P (u)\-$category\/} $\big( C_\F (u),C_{\hat\F^{^{\rm sc}}} (u)\big)\,;$ thus, as in~£2.1 above, {\it via\/}~the structural functor $C_{\P} (u)\to C_{\F} (u)$ we get a {\it regular central 
$k^*\-$extension\/} $C_{\hat\P^{^{\rm sc}}}(u)$ of  $C_{\P^{^{\rm sc}}} (u)$ and, in [10,~9.2.5], we have defined a {\it general decomposition map\/}
$$\delta^u_{\hat\P^{^{\rm sc}}} : \G_\K (\hat\P^{^{\rm sc}}) \too \G_k \big(C_{\hat\P^{^{\rm sc}}}(u)\big)
\eqno £3.9.1;$$
finally,  in [10,~Theorem~9.3] we state that the set of these {\it general decomposition maps\/} when $u$ runs over $\U$
 determines a $\K\-$module isomorphism
$$\K\otimes_\O \G_\K (\hat\P^{^{\rm sc}})\cong 
\prod_{u\in \U} \K\otimes_\O \G_k \big(C_{\hat\P^{^{\rm sc}}}(u)\big)
\eqno £3.9.2.$$

\medskip
\noindent
{\bf Theorem~£3.10.} {\it With the notation above, any local element $v$ of $\hat G(\hat\P^{^{\rm sc}})$
has a $\hat G(\hat\P^{^{\rm sc}})\-$conjugate in the image of $P$ and the restriction induces  a $\K\-$module isomorphism
$$\K\otimes_\O\G_k\big(C_{\hat G(\hat\P^{^{\rm sc}})} (v)\big) \cong 
\prod_{u\in \U_v} \K\otimes_\O \G_k \big(C_{\hat\P^{^{\rm sc}}}(u)\big)
\eqno £3.10.1\phantom{.}$$
where $\U_v$ denotes the set of $u\in \U$ such that the image is $\hat G(\hat\P^{^{\rm sc}})\-$conjugate to~$v\,.$\/}

\medskip
\noindent
{\bf Proof:} With  notation in~£1.9 and~£2.25, for any $u\in \U$  let us denote by $u^*$ the $p\-$element 
$\frak f_{\hat\P^{^{\rm sc}}}\big(\hat\tau_P (u)\big)$ in $\hat G(\hat\P^{^{\rm sc}})\,;$ it is quite clear that the 
$k^*\-$functor  $\frak f_{\hat\P^{^{\rm sc}}}\,\colon \hat\P^{^{\rm sc}}\to \hat G(\hat\P^{^{\rm sc}})$ still induces 
a $k^*\-$functor from $C_{\hat\P^{^{\rm sc}}}(u)$ to $C_{\hat G(\hat\P^{^{\rm sc}})} (u^*)\,;$  then, by restriction, we get  $\O\-$module
homomorphisms
$$\eqalign{F_u : \G_\K\big(C_{\hat G(\hat\P^{^{\rm sc}})} (u^*)\big) &\too 
\G\big(\K_*C_{\hat\P^{^{\rm sc}}}(u)\-\mod\big) \cr
f_u :\G_k\big(C_{\hat G(\hat\P^{^{\rm sc}})} (u^*)\big) &\too \G\big(k_*C_{\hat\P^{^{\rm sc}}}(u)\-\mod\big) \cr}
\eqno £3.10.2.$$

\smallskip
At this point, we claim that the following diagram is commutative 
$$\matrix{\G_\K \big(\hat G(\hat\P^{^{\rm sc}})\big)&\cong & \G (\K_*\hat\P^{^{\rm sc}}\-\mod)
&\cong &\G_\K (\hat\P^{^{\rm sc}})\cr
\hskip-40pt{\scriptstyle \delta^{u^*}_{\hat G(\hat\P^{^{\rm sc}})}}\hskip5pt\big\downarrow&&
&&\big\downarrow\hskip5pt{\scriptstyle \delta^{u^*}_{\hat\P^{^{\rm sc}}}}\hskip-25pt\cr
\G_k \big(C_{\hat G(\hat\P^{^{\rm sc}})}(u^*)\big)&\buildrel f_u\over\too& \G \big(k_*C_{\hat\P^{^{\rm sc}}}(u)\-\mod\big)
&\buildrel \cat_k^u\over\too &\G_k \big(C_{\hat\P^{^{\rm sc}}}(u)\big)\cr}
\eqno £3.10.3,$$
the isomorphisms in the top line coming from Corollaries~£2.23 and~£2.26. First of all,  since $C_{\hat\P^{^{\rm sc}}}(u)$ is a $k^*\-$subcategory of $\hat\P^{^{\rm sc}}\,,$ it is quite clear that the restriction determines the following commutative diagram
$$\matrix{\G_\K \big(\hat G(\hat\P^{^{\rm sc}})\big)&\cong & \G (\K_*\hat\P^{^{\rm sc}}\-\mod)\cr
\big\downarrow&&\big\downarrow\cr
\G_\K \big(C_{\hat G(\hat\P^{^{\rm sc}})}(u^*)\big)&\cong & \G (\K_*C_{\hat\P^{^{\rm sc}}}(u)\-\mod)\cr}
\eqno £3.10.4.$$
\eject

\smallskip
Moreover, it follows from [10,~4.6] that we also have restriction maps
$$\G_\K (\hat\P^{^{\rm sc}})\too \G_\K \big(C_{\hat\P^{^{\rm sc}}} (u)\big)\qq
\G_k (\hat\P^{^{\rm sc}})\too \G_k \big(C_{\hat\P^{^{\rm sc}}} (u)\big)
\eqno £3.10.5;$$
then, from the  definition of the left-hand restriction, it is not difficult to check that we also get a commutative diagram
$$\matrix{\G (\K_*\hat\P^{^{\rm sc}}\-\mod) &\buildrel \cat_\K\over\cong &\G_\K (\hat\P^{^{\rm sc}})\cr
\big\downarrow&&\big\downarrow\cr
\G \big(\K_*C_{\hat\P^{^{\rm sc}}}(u)\-\mod \big) &\buildrel \cat_\K^u\over\cong 
&\G_\K \big(C_{\hat\P^{^{\rm sc}}}(u)\big)\cr}
\eqno £3.10.6;$$
hence, we finally get the commutative diagram
$$\matrix{\G_\K \big(\hat G(\hat\P^{^{\rm sc}})\big)&\cong & \G (\K_*\hat\P^{^{\rm sc}}\-\mod)
&\buildrel \cat_\K\over\cong &\G_\K (\hat\P^{^{\rm sc}})\cr
\big\downarrow\hskip5pt {\scriptstyle {\rm Res}^{\hat G(\hat\P^{^{\rm sc}})}_{C_{\hat G(\hat\P^{^{\rm sc}})} (u^*)}}
\hskip-60pt&&\big\downarrow&&\big\downarrow\cr
\G_\K\big(C_{\hat G(\hat\P^{^{\rm sc}})} (u^*)\big) &\buildrel F_u\over\too 
&\G\big(\K_*C_{\hat\P^{^{\rm sc}}}(u)\-\mod\big)&\buildrel \cat_\K^u\over\cong 
&\G_\K \big(C_{\hat\P^{^{\rm sc}}}(u)\big)\cr}
\eqno £3.10.7.$$

\smallskip
On the other hand, as in~£3.6 above, it follows from the definition of the  {\it Brauer decomposition map\/} in [10,~3.4.2]
that we have the following commutative diagram
$$\matrix{\G_\K\big(C_{\hat G(\hat\P^{^{\rm sc}})} (u^*)\big) &\buildrel F_u\over\too 
&\G\big(\K_*C_{\hat\P^{^{\rm sc}}}(u)\-\mod\big)&\buildrel \cat_\K^u\over\cong 
&\G_\K \big(C_{\hat\P^{^{\rm sc}}}(u)\big)\cr
\big\downarrow\hskip5pt {\scriptstyle \delta_{C_{\hat G(\hat\P^{^{\rm sc}})} (u^*)}}\hskip-40pt&&\phantom{\Big\downarrow}&&\hskip-30pt{\scriptstyle \delta_{C_{\hat\P^{^{\rm sc}}} (u)}}\hskip5pt\big\downarrow\cr
\G_k\big(C_{\hat G(\hat\P^{^{\rm sc}})} (u^*)\big) &\buildrel f_u\over\too &\G\big(k_*C_{\hat\P^{^{\rm sc}}}(u)\-\mod\big) &\buildrel \cat_\K^u\over\too 
&\G_k \big(C_{\hat\P^{^{\rm sc}}}(u)\big)\cr}
\eqno £3.10.8.$$
But, the vertical left-hand arrow in diagram~£3.10.3 is the composition on the right- and on the left-hand of 
$\omega^{u^*}_{C_{\hat G(\hat\P^{^{\rm sc}})} (u^*)}$ with the left-hand arrows of diagrams~£3.10.7 and~£3.10.8
(cf.~£3.7.2); analogously,  the vertical right-hand arrow in diagram~£3.10.3 is the composition on the right- and on the left-hand 
of the following $\O\-$module automorphism defined in~[10,~9.2.4]
$$\Omega_{C_{\hat\P^{^{\rm sc}}} (u)}^u : \G_\K \big(C_{\hat\P^{^{\rm sc}}}(u)\big)\cong 
\G_\K \big(C_{\hat\P^{^{\rm sc}}}(u)\big)
\eqno £3.10.9\phantom{.}$$
with the right-hand arrows of diagrams~£3.10.7 and~£3.10.8.

\smallskip
Consequently, in order to show the commutativity of diagram~£3.10.3 it remains to prove the  commutativity of
the following diagram
$$\matrix{\G_\K\big(C_{\hat G(\hat\P^{^{\rm sc}})} (u^*)\big) &\buildrel F_u\over\too 
&\G\big(\K_*C_{\hat\P^{^{\rm sc}}}(u)\-\mod\big)&\buildrel \cat_\K^u\over\cong 
&\G_\K \big(C_{\hat\P^{^{\rm sc}}}(u)\big)\cr
\big\downarrow\hskip5pt {\scriptstyle \omega^{u^*}_{C_{\hat G(\hat\P^{^{\rm sc}})} (u^*)}}\hskip-40pt&&\phantom{\Big\downarrow}&&\hskip-30pt{\scriptstyle \Omega_{C_{\hat\P^{^{\rm sc}}} (u)}^u}\hskip5pt\big\downarrow\cr
\G_\K\big(C_{\hat G(\hat\P^{^{\rm sc}})} (u^*)\big) &\buildrel F_u\over\too 
&\G\big(\K_*C_{\hat\P^{^{\rm sc}}}(u)\-\mod\big)&\buildrel \cat_\K^u\over\cong 
&\G_\K \big(C_{\hat\P^{^{\rm sc}}}(u)\big)\cr}
\eqno £3.10.10;$$
for this purpose, recall that $\G_\K \big(C_{\hat\P^{^{\rm sc}}}(u)\big)$ is the {\it inverse limit\/} of the family
$\Big\{\G_\K \Big(\big(C_{\hat\P^{^{\rm sc}}}(u)\big)(\hat\frak q)\Big)\Big\}_{\hat\frak q}$ and then that the $\O\-$module automorphism $\Omega_{C_{\hat\P^{^{\rm sc}}} (u)}^u$ is the 
{\it inverse limit\/} of the family of $\O\-$module automorphisms
$$\omega_{(C_{\hat\P^{^{\rm sc}}}(u))(\hat\frak q)}^u : \G_\K \Big(\big(C_{\hat\P^{^{\rm sc}}}(u)\big)(\hat\frak q)\Big)\cong \G_\K \Big(\big(C_{\hat\P^{^{\rm sc}}}(u)\big)(\hat\frak q)\Big)
\eqno £3.10.11\phantom{.}$$
 where $\hat\frak q\,\colon \Delta_n\to C_{\hat\P^{^{\rm sc}}}(u)$ runs over all the 
 {\it $C_{\hat\P^{^{\rm sc}}}(u)\-$chains\/} (cf.~£2.1); that is to say, this family defines a {\it natural automorphism\/}
 [10,~9.2.3]
 $$\omega_{C_{\hat\P^{^{\rm sc}}}(u)}^u : \frak g_\K\circ \aut_{C_{\hat\P^{^{\rm sc}}}(u)} \cong 
 \frak g_\K\circ \aut_{C_{\hat\P^{^{\rm sc}}}(u)} 
 \eqno £3.10.12\phantom{.}$$
 and we set [10,~9.2.4]
 $$\Omega_{C_{\hat\P^{^{\rm sc}}} (u)}^u = \lim_{\longleftarrow}\,(\omega_{C_{\hat\P^{^{\rm sc}}}(u)}^u )
 \eqno £3.10.13.$$
 
 \smallskip
 But, for any {\it $C_{\hat\P^{^{\rm sc}}}(u)\-$chain\/} $\hat\frak q\,\colon \Delta_n\to C_{\hat\P^{^{\rm sc}}}(u)$
 we have the structural $\O\-$module homomorphism 
 $$\ell^u_{\hat\frak q} : \G_\K \big(C_{\hat\P^{^{\rm sc}}}(u)\big)\too
 \G_\K \Big(\big(C_{\hat\P^{^{\rm sc}}}(u)\big)(\hat\frak q)\Big)
 \eqno £3.10.14\phantom{.}$$
  and, moreover, the $k^*\-$functor from $C_{\hat\P^{^{\rm sc}}}(u)$ to $C_{\hat G(\hat\P^{^{\rm sc}})} (u^*)$ induced by 
$\frak f_{\hat\P^{^{\rm sc}}}$ still induces a $k^*\-$group homomorphism 
$$\theta^u_{\hat \frak q} : \big(C_{\hat\P^{^{\rm sc}}}(u)\big)(\hat\frak q)\too C_{\hat G(\hat\P^{^{\rm sc}})} (u^*)
\eqno £3.10.15;$$
then, it is easily checked that we have the following commutative diagram
$$\matrix{\G_\K\big(C_{\hat G(\hat\P^{^{\rm sc}})} (u^*)\big) \buildrel F_u\over\too 
&\G\big(\K_*C_{\hat\P^{^{\rm sc}}}(u)\-\mod\big)&\buildrel \cat_\K^u\over\cong 
\G_\K \big(C_{\hat\P^{^{\rm sc}}}(u)\big)\cr
{\atop {\rm Res}_{\theta^u_{\hat \frak q}}}\searrow\hskip-50pt&
\phantom{\Big\downarrow}&\hskip-50pt\swarrow{\atop \ell^u_{\hat\frak q}}\cr
&\G_\K \Big(\big(C_{\hat\P^{^{\rm sc}}}(u)\big)(\hat\frak q)\Big)\cr}
\eqno £3.10.16.$$
Finally, the commutativity of diagram~£3.10.10 follows from definition~£3.10.13 and the obvious commutativity of the following 
diagram
$$\matrix{\G_\K\big(C_{\hat G(\hat\P^{^{\rm sc}})} (u^*)\big) &\buildrel {\rm Res}_{\theta^u_{\hat \frak q}}\over\too 
&\G_\K \Big(\big(C_{\hat\P^{^{\rm sc}}}(u)\big)(\hat\frak q)\Big)\cr
\hskip-50pt{\scriptstyle \omega^{u^*}_{C_{\hat G(\hat\P^{^{\rm sc}})} (u^*)}}\hskip5pt\big\downarrow&
&\big\downarrow\hskip5pt {\scriptstyle \omega_{(C_{\hat\P^{^{\rm sc}}}(u))(\hat\frak q)}^u}\hskip-50pt \cr
\G_\K\big(C_{\hat G(\hat\P^{^{\rm sc}})} (u^*)\big) &\buildrel {\rm Res}_{\theta^u_{\hat \frak q}}\over\too 
&\G_\K \Big(\big(C_{\hat\P^{^{\rm sc}}}(u)\big)(\hat\frak q)\Big)\cr}
\eqno £3.10.17.$$

\smallskip
Now, the commutativity of diagram~£3.10.3 for any $u\in \U$ determines the following commutative diagram (cf.~£3.9.2)
$$\matrix{\K\otimes_\O\G_\K \big(\hat G(\hat\P^{^{\rm sc}})\big)&\cong &\K\otimes_\O \G_\K (\hat\P^{^{\rm sc}})\cr
\big\downarrow&\phantom{\Big\downarrow}&\wr\Vert\cr
{\displaystyle \prod_{u\in \U}}\K\otimes_\O \G_k\big(C_{\hat G(\hat\P^{^{\rm sc}})} (u^*)\big)
&\too & {\displaystyle \prod_{u\in \U}} \K\otimes_\O \G_k \big(C_{\hat\P^{^{\rm sc}}}(u)\big)\cr}
\eqno £3.10.18$$
\eject
\noindent
which forces the vertical left-hand arrow to be injective and the bottom arrow to be surjective; then, according to 
isomorphism~£3.8.2, the injectivity of the vertical left-hand arrow implies that the set $\{u^*\}_{u\in \U}$ contains 
a representative for any conjugacy class of {\it local elements\/} of $\hat G(\hat\P^{^{\rm sc}})\,;$ moreover, the surjectivity 
of the bottom arrow implies that the restriction to the image of the  vertical left-hand arrow is an isomorphism
and then the announced isomorphism~£3.10.1 follows easily.

\bigskip
\noindent
{\bf Corollary~£3.11.}{\it The $\O\-$module homomorphism $\cat_k \colon \G (k_*\hat\P^{^{\rm sc}}\!\!\!\-\mod)
\to \G_k (\hat\P^{^{\rm sc}})$ is bijective.\/}

\medskip
\noindent
{\bf Proof:} Since the  the vertical homomorphisms in  diagram~£3.6.2 above are surjective, this homomorphism is surjective.
Moreover, it follows from isomorphism~£3.10.1 for $v= 1$ that we have the  $\K\-$module isomorphism
$$\K\otimes_\O\G_k\big(\hat G(\hat\P^{^{\rm sc}})\big) \cong \K\otimes_\O \G_k (\hat\P^{^{\rm sc}})
\eqno £3.11.1\phantom{.}$$
and therefore, according to the $\O\-$module isomorphism~£3.5.2, $\cat_k$ is also injective.

 \bigskip
\noindent
{\bf References}
\medskip
\noindent
{\cds[1]  Richard Brauer, {\cdt On blocks and sections in finite groups, I and II\/}, Amer. J. Math. 89(1967), 1115-1136,
90(1968), 895-925
\smallskip
\noindent
[2] Michel Brou\'e, {\cdt Radical, hauteurs, p-sections et blocs\/}, Ann. of Math. 107(1978), 89-107
\smallskip\noindent
[3] Andrew Chermak. {\cdt Fusion systems and localities\/}, 
Acta Mathematica, 211(2013), 47-139.
\smallskip
\noindent
[4] Daniel Gorenstein, {\cdt ``Finite Groups''\/}, Harper's Series, 1968, Harper and Row
\smallskip\noindent
[5] Bob Oliver. {\cdt Existence and Uniqueness of Linking Systems: Chermak's proof via obstruction theory\/}, 
Acta Mathematica, 211(2013), 141-175.
\smallskip\noindent
[6]\phantom{.} Llu\'\i s Puig, {\cdt Local fusions in block source algebras\/},
Journal of Algebra, 104(1986), 358-369. 
\smallskip
\noindent
[7] Llu\'\i s Puig, {\cdt Pointed groups and  construction of modules}, Journal of Algebra, 116(1988), 7-129
\smallskip\noindent
[8]\phantom{.}  Llu\'\i s Puig, {\cdt Nilpotent blocks and their source
algebras}, Inventiones math., 93(1988), 77-116.
\smallskip
\noindent
[9] Llu\'\i s Puig, {\cdt ``Frobenius Categories versus Brauer Blocks''\/}, Progress in Math. 274, 2009, Birkh\"auser, Basel
\smallskip
\noindent
[10] Llu\'\i s Puig, {\cdt Ordinary Grothendieck groups of a Frobenius P-category\/}, Algebra Colloquium 18(2011), 1-76
\smallskip
\noindent
[11] Llu\'\i s Puig, {\cdt Existence, uniqueness and functoriality of the perfect
locality over a Frobenius P-category\/}, arxiv.org/abs/1207.0066, submitted to Algebra Colloquium.
\smallskip\noindent
[12] Llu\'\i s Puig, {\cdt A criterion on trivial homotopy\/},  arxiv.org/abs/1308.3765, submitted to Journal of Algebra.
\smallskip
\noindent
[13] Llu\'\i s Puig, {\cdt Affirmative answer to a question of Linckelmann\/}, arxiv.org/abs/1507.04278, submitted to 
Journal of Algebra.}

\end